\documentclass{amsart}
\usepackage{amsmath}%
\usepackage{amsfonts}%
\usepackage{amssymb}%
\usepackage{graphicx}

\usepackage{hyperref}
\usepackage{enumerate}
\usepackage{eucal}

\newtheorem{theorem}{Theorem}
\newtheorem*{theorem*}{Theorem}
\theoremstyle{plain}

\newtheorem{claim}{Claim}

\newtheorem{conjecture}{Conjecture}
\newtheorem{corollary}{Corollary}

\newtheorem{definition}{Definition}
\newtheorem{example}{Example}

\newtheorem{lemma}{Lemma}
\newtheorem*{lemma*}{Lemma}

\newtheorem{proposition}{Proposition}
\newtheorem{remark}{Remark}

\numberwithin{equation}{section}

\newtheorem{remarks}{Remarks}

\newcommand{\R}{\mathbb{R}}
\newcommand{\Sp}{\mathbb{S}}
\newcommand{\Hy}{\mathbb{H}}
\newcommand{\Q}{\mathbb{Q}}

\newcommand{\Span}{\mbox{\normalfont{span}}}

\newcommand{\trac}{\mbox{\normalfont{tr\,}}}

\newcommand{\beeq}{\begin{equation}}
\newcommand{\eneq}{\end{equation}}
\newcommand{\beeqs}{\begin{eqnarray*}}
\newcommand{\eneqs}{\end{eqnarray*}}
\newcommand{\besp}{\begin{split}}
\newcommand{\ensp}{\end{split}}
\newcommand{\bepr}{\begin{proof}}
\newcommand{\enpr}{\end{proof}}
\newcommand{\bethr}{\begin{theorem}}
\newcommand{\enthr}{\end{theorem}}
\newcommand{\beths}{\begin{theorem*}}
\newcommand{\enths}{\end{theorem*}}
\newcommand{\becor}{\begin{corollary}}
\newcommand{\encor}{\end{corollary}}
\newcommand{\bere}{\begin{remark}}
\newcommand{\enre}{\end{remark}}
\newcommand{\bers}{\begin{remarks}}
\newcommand{\enrs}{\end{remarks}}
\newcommand{\beres}{\begin{remark*}}
\newcommand{\enres}{\end{remark*}}
\newcommand{\bele}{\begin{lemma}}
\newcommand{\enle}{\end{lemma}}
\newcommand{\beles}{\begin{lemma*}}
\newcommand{\enles}{\end{lemma*}}
\newcommand{\bepro}{\begin{proposition}}
\newcommand{\enpro}{\end{proposition}}
\newcommand{\bepros}{\begin{proposition*}}
\newcommand{\enpros}{\end{proposition*}}
\newcommand{\becl}{\begin{claim}}
\newcommand{\encl}{\end{claim}}
\newcommand{\beex}{\begin{example}}
\newcommand{\enex}{\end{example}}
\newcommand{\beexs}{\begin{example*}}
\newcommand{\enexs}{\end{example*}}
\newcommand{\beco}{\begin{conjecture}}
\newcommand{\enco}{\end{conjecture}}
\newcommand{\becos}{\begin{conjecture*}}
\newcommand{\encos}{\end{conjecture*}}
\newcommand{\bede}{\begin{definition}}
\newcommand{\bedes}{\begin{definition*}}

\newcommand{\lp}{\left(}
\newcommand{\rp}{\right)}
\newcommand{\lrb}{\left[}
\newcommand{\rrb}{\right]}
\newcommand{\lb}{\left\{}
\newcommand{\rb}{\right\}}
\newcommand{\la}{\left\langle}
\newcommand{\ra}{\right\rangle}

\textheight=\dimexpr
  \ifcase 10pt
    \or 
    \or 
    \or 
    \or 
    \or 
    \or 
    \or 
    \or 
    57\or 
    52\or 
    48\or 
    44\or 
    41\fi 
  \baselineskip+\topskip\relax
\begin{document}
\title[SPHERICAL 2-DUPIN SUBMANIFOLDS]{SPHERICAL 2-DUPIN SUBMANIFOLDS}
\author[Antonio J. Di Scala]{Antonio J. Di Scala}
\address{Antonio J. Di Scala -- Politecnico di Torino\newline%
\indent Corso Duca degli Abruzzi, 24, 10129 Turin, Italy}%
\email{antonio.discala@polito.it}%

\author[Guilherme Machado de Freitas]{Guilherme Machado de Freitas}
\address{Guilherme Machado de Freitas -- Politecnico di Torino\newline%
\indent Corso Duca degli Abruzzi, 24, 10129 Turin, Italy}%
\email{guimdf1987@icloud.com}%

\thanks{A.J. Di Scala is member of GNSAGA of INdAM}
\thanks{G. Machado de Freitas's research was partially supported by CNPq/Brazil}

\subjclass[2010]{Primary 53B25; Secondary 53C29} %
\keywords{2-Dupin submanifolds, standard embeddings of $\mathbb{FP}^2$, Lie sphere geometry}%

\begin{abstract}
We show that every spherical 2-Dupin submanifold that is not a hypersurface is conformally congruent to the standard embedding of the real, complex, quaternionic or octonionic projective plane. We also classify 2-CPC, 2-umbilical and weakly 2-umbilical submanifolds in space forms.
\end{abstract}

\maketitle

\section{Introduction}
An isometric immersion $f:M^n\to\tilde{M}^m$ is said to be \emph{umbilical} at $x\in M^n$ if there exists $\eta$ in its normal space $N_fM\lp x\rp$ such that
its second fundamental form
\beeqs
\alpha\lp X,Y\rp=\la X,Y\ra\eta
\eneqs
for all $X,Y\in T_xM$. Clearly, in this case $\eta$ is the mean curvature vector $\mathcal{H}\lp x\rp$ of $f$ at $x$. Equivalently, $f$ is umbilical at $x$ if its shape operator
\beeq\label{usso}
A_\xi=\la\mathcal{H}\lp x\rp,\xi\ra I
\eneq
for every $\xi\in N_fM\lp x\rp$. A submanifold is called \emph{umbilical} if it is umbilical at every point.

An isometric immersion $f:M^n\to\tilde{M}^m$ is said to have \emph{parallel mean curvature vector field} if
\beeq\label{pmcvf}
\nabla^\perp_X\mathcal{H}=0
\eneq
for all $x\in M^n$ and $X\in T_xM$.

In particular, if $f$ has parallel mean curvature vector field then it follows that $\left\|\mathcal{H}\right\|$ is constant along $M^n$.

An immersion is called an \emph{extrinsic sphere} if it is umbilical and the mean curvature vector field is parallel.

It is a well-known fact that every $n$-dimensional (complete) extrinsic sphere of the space form $\Q^m_{\tilde{c}}$ is isometric to $\Q^n_c$, where $c=\tilde{c}+\left\|\mathcal{H}\right\|^2$, and it is always a hypersurface in the sense that its codimension can be reduced to one, i.e., it is contained in a totally geodesic submanifold $\Q^{n+1}_{\tilde{c}}$. Moreover, it is uniquely determined by and can be explicitly described in terms of one of its points together with its corresponding tangent space and mean curvature vector. Indeed, in Euclidean space, for instance, apart from affine subspaces (totally geodesic submanifolds), any extrinsic sphere is an actual round sphere, while in the spherical and hyperbolic space forms (regarded as round spheres of the Euclidean and Lorentzian space, respectively) an extrinsic sphere is always the intersection with an affine subspace of the ambient Euclidean or Lorentzian space, respectively. Of course an extrinsic sphere is totally geodesic if and only if $\mathcal{H}=0$. We refer to Dajczer \cite{dajczer1990submanifolds} or Berndt-Console-Olmos \cite{MR3468790} for the many details of this result.

It is also well known that the parallelism \eqref{pmcvf} of the mean curvature vector field is automatically satisfied for $n$-dimensional umbilical submanifolds in space forms with $n\geq2$. This is a simple consequence of the Codazzi equation. In other words, any $n$-dimensional umbilical submanifold of a space form is always an extrinsic sphere if $n\geq2$. On the other hand, this is clearly not true if $n=1$, since every curve is trivially umbilical.

We now discuss how the definitions of umbilical submanifolds and extrinsic spheres can be generalized to encompass some of the most beautiful and widely studied classes of submanifolds, but before proceeding we agree that manifolds are always connected throughout the article and we denote by
\beeqs
UN_f^{m-1}=\lb\lp x,\xi\rp\in N_fM:\left\|\xi\right\|=1\rb
\eneqs
the unit normal bundle of a given isometric immersion $f:M^n\to\tilde{M}^m$, except when $f$ is an oriented hypersurface, in which case the right-hand side of the above equation is not connected and so $UN_f^n$ will stand for one of its two connected components. We also use the natural identification of $T_{\lp x,\xi\rp}UN_f$ as a subspace of $T_xM\times T_{f\lp x\rp}\tilde{M}$ by means of differentiation with respect to the Levi-Civita connection $\tilde{\nabla}$ of $\tilde{M}$, i.e., given a curve $\lp x\lp t\rp,\xi\lp t\rp\rp$ in $UN_f^{m-1}$, we identify
\beeq\label{wfeq}
\lp x,\xi\rp'\simeq\lp x',\tilde{\nabla}_{x'}\xi\rp=\lp x',-f_*A_\xi x'+\nabla^\perp_{x'}\xi\rp,
\eneq
where we have used the Weingarten formula in the equality. Since $\xi\lp t\rp$ is a unit vector field along $x\lp t\rp$, it follows that $\nabla^\perp_{x'}\xi$ is orthogonal to $\xi$, and hence
\beeqs
T_{\lp x,\xi\rp}UN_f\simeq\lb\lp X,-f_*A_\xi X+\eta\rp:X\in T_xM\text{ and }\eta\in N_fM\lp x\rp,\eta\cdot\xi=0\rb.
\eneqs

Note that equation \eqref{usso} can be simply translated into the fact that $f$ has exactly one principal curvature with respect to any $\lp x,\xi\rp\in UN_f^{m-1}$. It is hence quite natural to extend this definition to the case in which $f$ has precisely $k$ principal curvatures with respect to every $\lp x,\xi\rp\in UN_f^{m-1}$. An isometric immersion $f:M^n\to\tilde{M}^m$ is said to be $k$\emph{-umbilical} at $x\in M^n$ if it has exactly $k$ distinct principal curvatures with respect to all $\lp x,\xi\rp\in UN_f^{m-1}$. A submanifold is called $k$\emph{-umbilical} if it is $k$-umbilical at every point. In this case, each principal curvature $\kappa_i:UN_f^{m-1}\to\R,\ 1\leq i\leq k,$ has constant multiplicity and is then smooth as well as its corresponding principal distribution
\beeqs
E_i\lp x,\xi\rp=\lb X\in T_xM:A_\xi X=\kappa_iX\rb,
\eneqs
by the general theory of symmetric smooth tensors on Riemannian manifolds (see, for example, Cecil-Ryan \cite{cecil1978}, Nomizu \cite{doi:10.1080/03081087308817014} or Ryan \cite{ryan1969}). Of course, $E_i$ is not a distribution in the strict sense of the word, although we will soon associate it to its horizontal lift $T_i$, which is then literally a distribution on $UN_f^{m-1}$. Notice also that, if $k\geq2$ and $f:M^n\to\Q^{n+p}_c$ is $k$-umbilical at some point $x$, the definition forces $f$ to be conformally substantial, since in this case no umbilical normal direction at $x$ will then exist.

The definition of an extrinsic sphere, on the other hand, is a little bit more complex, since it involves not only the algebraic nature of the second fundamental form, but also the differential requirement that the mean curvature vector field be parallel. The attempt to generalize the latter condition to $k$-umbilical submanifolds, depending on how we look at it, will lead us through a nested sequence
\beeq\label{nsc}
\mathfrak{U}_k\subset\mathfrak{I}_k\subset\mathfrak{D}_k
\eneq
of three remarkable classes of submanifolds.

Let $S^n$ be an extrinsic sphere in a space form. As discussed at the beginning of the section, we can assume that $S^n$ is an orientable hypersurface. Then the whole information of its second fundamental form is contained in its shape operator $A$ with respect to one of its two unit normal vector fields and condition \eqref{pmcvf} reduces to the constancy of its unique principal curvature (which is no longer than its mean curvature).

The simplest and most obvious way of extending the last property to a $k$-umbilical submanifold is by requiring that it have the same principal curvatures with respect to all normal directions. A $k$-umbilical isometric immersion $f:M^n\to\tilde{M}^m$ is said to be $k$\emph{-unipotent} if its principal curvatures $\kappa_i:UN_f^{m-1}\to\R,\ 1\leq i\leq k,$ are constant. This is a very special symmetry condition that imposes strong geometric restrictions on $f$ even punctually, since it will then shape up the same way in all directions. Prototypical examples of unipotent submanifolds of a Riemannian manifold $\tilde{M}^m$ are given by the so-called normal homogeneous submanifolds, that is, orbits of a connected closed subgroup of $I^0\lp\tilde{M}\rp$ whose slice representation at a point acts transitively on the unit sphere in the normal space at that point.
\bere
In \cite{DiScala2006}, the first author used the notion of unipotency in the sense that the shape operator $A_\xi$ satisfies the equation $A_\xi^2=k\left\|\xi\right\|^2I$, where $k>0$ is constant, which except for hypersurfaces is easily seen to be equivalent to our definition of a 2-unipotent submanifold.
\enre
The second possibility is to require that the principal curvatures be constant just along parallel normal vector fields. Since these vector fields do not exist even locally in general, we must pick them always along piecewise differentiable curves, and this brings us to the classical concepts of submanifolds with constant principal curvatures and isoparametric submanifolds. Recall that an isometric immersion $f:M^n\to\tilde{M}^m$ is said to have constant principal curvatures if for any parallel normal vector field $\xi\lp t\rp$, along any piecewise differentiable curve, the principal curvatures in direction $\xi\lp t\rp$ are constant. If in addition the normal bundle of $f$ is flat, one says that the submanifold is isoparametric.
\bere
Isoparametric hypersurfaces were studied by many prominent mathematicians in the 20th century, including B. Segre, Levi-Civita and \'{E}. Cartan, and they are still an active research field. The beginning of the history of their generalizations to higher codimensions goes back to the 80s, with the (sometimes independent) work of many authors: J. Eells \cite{MR825258}, D. Gromoll and K. Grove \cite{Gromoll1985}, Q. M. Wang \cite{MR1024799}, C. E. Harle \cite{Harle1982}, W. Str\"{u}bing \cite{Strübing1986} and C.-L. Terng \cite{terng1985}. There are different aspects of these generalizations that are actually strictly related: isoparametric maps, isoparametric submanifolds and submanifolds with constant principal curvatures (for more details on the historical development see \cite{MR1736861}). The general notion of an isoparametric map is credited to C.-L. Terng. It turns out that regular level submanifolds of an isoparametric map are isoparametric submanifolds. Conversely, any isoparametric submanifold $M^n$ of $\R^m$ determines a polynomial isoparametric map on $\R^m$, which has $M^n$ as a regular level set (cf. \cite{MR972503} and \cite{MR3468790}). The notion of isoparametric submanifold is nowadays also regarded as originally given by C.-L. Terng, even though it was first stated by C. E. Harle in \cite{Harle1982}. There is a strong link between those submanifolds and $s$-representations. In fact, as a consequence of a result by Thorbergsson \cite{10.2307/2944343}, the orbits of $s$-representations are almost all submanifolds with constant principal curvatures. It was conjectured by C. Olmos \cite{olmos1994}, twenty years ago, that a full and irreducible homogeneous submanifold of the sphere such that the normal holonomy group is not transitive must be an orbit of an irreducible $s$-representation. The conjecture was already known to be true for surfaces, and recently C. Olmos and R. Ria\~{n}o-Ria\~{n}o \cite{olmos2015} have proved it also in dimension three. For many other interesting facts and results on isoparametric and homogeneous submanifolds, the interested reader can refer to Heintze-Olmos-Thorbergsson \cite{doi:10.1142/S0129167X91000107} and Berndt-Console-Olmos \cite{MR3468790}.
\enre
Since the condition of constant principal curvatures is far more important than having flat normal bundle and actually captures the essence of an isoparametric submanifold, we will make no distinction between them and submanifolds with constant principal curvatures, and both will be called CPC submanifolds for simplicity, whether their normal bundle is flat or not. It is also useful to indicate explicitly the constant number $k$ of principal curvatures of $f$ by saying that $f$ is a $k$\emph{-CPC submanifold}. Observe that, if $\lp x\lp t\rp,\xi\lp t\rp\rp\in UN_f^{m-1}$ is a parallel normal vector field, then equation \eqref{wfeq} yields
\beeqs
\lp x,\xi\rp'=\lp x',-f_*A_\xi x'\rp,
\eneqs
and thus the constancy of the principal curvatures $\kappa_i,\ 1\leq i\leq k,$ along $\lp x,\xi\rp$ is equivalent to
\beeqs
{\kappa_i}_*\lp x',-f_*A_\xi x'\rp=0.
\eneqs
In other words, $f$ is $k$-CPC if and only if every $\kappa_i$ is constant along the horizontal lift
\beeqs
T_{\lp x,\xi\rp}M=\lb\lp X,-f_*A_\xi X\rp:X\in T_xM\rb
\eneqs
of the tangent bundle.
\bere
Unlike our definition of a CPC submanifold, those found in the literature do not require the submanifold to be $k$-umbilical.
\enre
Finally, a last possible extension of the concept of an extrinsic sphere to $k$-umbilical submanifolds is obtained by requiring the even weaker condition that each principal curvature $\kappa_i,\ 1\leq i\leq k,$ be constant but only along its lines of curvature, i.e., parallel normal vector fields $\lp x\lp t\rp,\xi\lp t\rp\rp\in UN_f^{m-1}$ such that
\beeqs
A_\xi x'=\kappa_i x',
\eneqs
or equivalently, that each $\kappa_i$ be constant along the horizontal lift
\beeqs
T_i=\lb\lp X,-f_*A_\xi X\rp:X\in E_i\rb
\eneqs
of its principal distribution $E_i$. In codimension one, $T_i$ can be obviously identified with $E_i$ and the above notion reduces to the definition of a Dupin hypersurface. It is hence quite natural to call a general $k$-umbilical isometric immersion $f:M^n\to\tilde{M}^m$ a $k$\emph{-Dupin submanifold} if $f$ satisfies the above property.
\bere
See Reckziegel \cite{Reckziegel1979} for a general Riemannian treatment of the notion of principal curvatures and curvature submanifolds in the case of an isometric immersion $f:M^n\to\Q^{n+p}_c$ of codimension greater than one. In that case, Reckziegel defines a curvature submanifold to be a connected submanifold $S\subset M^n$ for which there is a parallel (with respect to the normal connection) section of the unit normal bundle $\eta:S\to UN_f^{n+p-1}$ such that for each $x\in S$, the tangent space $T_xS$ is equal to some eigenspace of $A_{\eta\lp x\rp}$. The corresponding principal curvature $\kappa:S\to\R$ is then a smooth function on $S$. Pinkall \cite{Pinkall1986} calls an isometric immersion $f:M^n\to\Q^{n+p}_c$ of codimension greater than one ``Dupin'' if all of its principal curvatures are Dupin, i.e., constant along each corresponding curvature submanifold (in the sense of Reckziegel). Pinkall's definition is clearly equivalent to our definition of a $k$-Dupin submanifold in the case where $f$ is $k$-umbilical. Furthermore, the notion of Dupin can be generalized even to the larger class of Legendre submanifolds, so that our definition of a Dupin submanifold can be expressed in terms of its Legendre lift, as we will see in Subsection \ref{lsgds}. The basic idea is that a submanifold of codimension greater than one is $k$-Dupin if and only if its Legendre lift is $\lp k+1\rp$-Dupin in the sense of Pinkall \cite{Pinkall1985}. In any case, we are surprised that the concept of a Dupin submanifold of higher codimension has not yet been much explored in the literature. On the other hand, we mention, for the record, the work \cite{dajczer2005} by Dajczer, Florit and Tojeiro, in which they studied a weak notion of reducibility for submanifolds that carry a so-called Dupin principal normal.
\enre
The link with Lie sphere geometry provided by Legendre lifts is at the core of our approach throughout this article and allows us to deduce several interesting properties of Dupin submanifolds from the corresponding ones satisfied by their lifts. It will follow for instance that Dupin submanifolds are invariant under conformal transformations of and between ambient spaces, and also that the constancy of a principal curvature along the leaves of its principal foliation is automatic when its multiplicity is greater than one.
\bere
Unlike our definition of a Dupin submanifold, those found in the literature for hypersurfaces also do not require the submanifold to be $k$-umbilical, and the term ``proper'' is reserved for when this is the case.
\enre
It is clear from the above definitions that every $k$-unipotent submanifold is $k$-CPC and that every $k$-CPC submanifold is in turn $k$-Dupin, so that \eqref{nsc} holds in the obvious notation. It is also clear that for $k=1$ the first notion only makes sense for hypersurfaces, in which case all three definitions collapse into that of an extrinsic sphere. Actually, in the case of hypersurfaces the first two concepts of $k$-unipotent and $k$-CPC hypersurfaces always coincide, for all $k\in\mathbb{N}$.

In the next case $k=2$, the class of 2-unipotent submanifolds in space forms has been studied by the first author in \cite{DiScala2006}, who has shown that in codimension greater than one it is nonempty essentially only in the spheres $\Sp^4,\Sp^7,\Sp^{13},\Sp^{25}$ and in each case consists of the sole standard embedding of the real, complex, quaternionic or octonionic projective plane, respectively. The reason for the `essentially' qualifier is that such embeddings can of course also be composed with the inverse of the \emph{similarity} $\theta_c:\Sp^m_c\to\Sp^m$ given by
\beeqs
\theta_c\lp x\rp=\sqrt{c}\,x,
\eneqs
where $\Sp^m_c\subset\R^{m+1}$ is regarded as
\beeqs
\Sp^m_c=\lb x\in\R^{m+1}:\left\|x\right\|^2=\frac{1}{c}\rb,
\eneqs
in order to yield the versions of the standard embeddings in $\Sp^4_c,\Sp^7_c,\Sp^{13}_c,\Sp^{25}_c$ instead.
\bethr[Di Scala \cite{DiScala2006}]\label{dsth}
Let $f:M^n\to\Q^{n+p}_c,\ p\geq2,$ be an isometric immersion. Then, the following facts are equivalent:
\begin{enumerate}[(i)]
\item $f$ is 2-unipotent;
\item $c>0,\ n=2^\mu,\ p=2^{\mu-1}+1$ for $\mu=1,2,3,4$ and, up to congruence and the similarity $\theta_c$, $f\lp M\rp$ is an open subset of $\psi_{\mathbb{F}}\lp\mathbb{FP}^2\rp$, where $\psi_{\mathbb{F}}:\mathbb{FP}^2\to\Sp^{3\cdot2^{\mu-1}+1}$ denotes the standard embedding of the projective plane $\mathbb{FP}^2$, for $\mathbb{F}=\R,\mathbb{C},\mathbb{H},\mathbb{O}$, according to whether $\mu=1,2,3,4$, respectively.
\end{enumerate}
\enthr
The standard embeddings of the projective planes are also familiar to us from Cartan's \cite{Cartan1939} classification of 3-isoparametric hypersurfaces, which together with the above remark on Legendre lifts suggests that they might exhaust the class of 2-CPC submanifolds in space forms too in codimension greater than one. This is actually the case, and is the content of one of our main results in this article.
\bethr\label{2ith}
Let $f:M^n\to\Q^{n+p}_c,\ p\geq2,$ be a connected 2-CPC isometric immersion. Then $c>0,\ n=2^\mu,\ p=2^{\mu-1}+1$ for $\mu=1,2,3,4$ and there exists an isometry $\phi$ of $\Sp^{3\cdot2^{\mu-1}+1}_c$ such that $\theta_c\circ\phi\circ f\lp M\rp$ is an open subset of the image of the standard embedding $\psi_{\mathbb{F}}:\mathbb{FP}^2\to\Sp^{3\cdot2^{\mu-1}+1}$ of the projective plane $\mathbb{FP}^2$, where $\mathbb{F}=\R,\mathbb{C},\mathbb{H},\mathbb{O}$, respectively, and $\theta_c:\Sp^{3\cdot2^{\mu-1}+1}_c\to\Sp^{3\cdot2^{\mu-1}+1}$ is the similarity.

Conversely, for every isometry $\phi$ of $\Sp^{3\cdot2^{\mu-1}+1}_c$, the composition $\phi\circ\theta_c^{-1}\circ\psi_{\mathbb{F}}:\mathbb{FP}^2\to\Sp^{3\cdot2^{\mu-1}+1}_c$ is a 2-CPC isometric immersion.
\enthr
\bere
The remaining case $p=1$ in Theorems \ref{dsth} and \ref{2ith} corresponds to the well-known classification of 2-isoparametric hypersurfaces in space forms. For $\R^{n+1}$ this classification is due to Levi-Civita \cite{levi1937famiglie} for $n=2$ and to Segre \cite{segre1938famiglie} in general, while for $\Hy^{n+1}$ and $\Sp^{n+1}$ it is due to E. Cartan \cite{Cartan1938}, \cite{Cartan1939}.

A 2-isoparametric hypersurface of $\R^{n+1}$ or $\Hy^{n+1}$ is always a spherical cylinder (or tube around a totally geodesic submanifold, depending on the reader's philosophical bent), whereas in $\Sp^{n+1}$ it is a Riemannian product of two spheres.
\enre
\bere
For another interesting characterization of the Veronese submanifolds $\psi_{\mathbb{F}}$ see the recent paper \cite{olmos2015} by C. Olmos and R. Ria\~{n}o-Ria\~{n}o. These submanifolds also appear implicitly or explicitly in several other contexts and works, such as in Console-Olmos \cite{Console1998} , Nurowski \cite{Nurowski20081148}, \cite{Nurowski:939061} and also \cite{MR3451462} by the second author, just to give a few references. The importance of the Veronese embeddings in the geometry of submanifolds cannot be sufficiently emphasized. It seems that time and again, when people least expect it, new important developments come out involving these simple and compellingly beautiful submanifolds.
\enre
Finally, concerning the last class of 2-Dupin submanifolds, in view of conformal invariance, many other spherical examples of these are obtained after a conformal deformation of one of the standard embeddings $\psi_{\mathbb{F}}$ and even in Euclidean and hyperbolic space forms via the stereographic projections $\pi_0:\R^m\to\Sp^m$ and $\pi_c:\Hy^m_c\to\Sp^m$ with pole $\sigma=\lp-1,0,\dots,0\rp$ given by
\beeq\label{sgpr}
\begin{gathered}
\pi_0\lp x\rp=\lp\frac{1-x\cdot x}{1+x\cdot x},\frac{2x}{1+x\cdot x}\rp,\\
\pi_c\lp x\rp=\lp\frac{1-y\cdot y}{1+y\cdot y},\frac{2y}{1+y\cdot y}\rp,\quad y=\lp x_2,\dots,x_{m+1}\rp/\lp x_1+c^{-1/2}\rp,
\end{gathered}
\eneq
where we regard $\Hy^m_c$ in the Lorentzian space $\mathbb{L}^{m+1}$ as
\beeqs
\Hy^m_c=\lb x\in\mathbb{L}^{m+1}:\left\|x\right\|^2=\frac{1}{c}\rb.
\eneqs
Nevertheless, the main goal of this paper is to prove that these are the only possible new examples of 2-Dupin submanifolds in space forms. Since these submanifolds are invariant under conformal transformations between ambient spaces, the statement is essentially the same whether it is considered in the sphere, in Euclidean space or in hyperbolic space. Among these three geometries, the spherical one is the best suited for displaying solutions of problems in conformal geometry with some Riemannian-geometric flavor, especially for the fact that the sphere has this distinguished property that any two given extrinsic spheres of the same dimension are always similar, while in the other two space forms also stereographic projections must be taken into account in order to move from one type of extrinsic sphere to another. Hence, we confine ourselves here to stating our classification of 2-Dupin submanifolds just in the unit sphere $\Sp^{m}$.
\bethr\label{omth}
Let $f:M^n\to\Sp^{n+p}$, $p\geq2,$ be a connected 2-Dupin isometric immersion. Then $n=2^\mu$, $p=2^{\mu-1}+1$ for $\mu=1,2,3,4$ and there exists a conformal transformation $\phi$ of $\Sp^{3\cdot2^{\mu-1}+1}$ such that $\phi\circ f\lp M\rp$ is an open subset of the image of the standard embedding $\psi_{\mathbb{F}}:\mathbb{FP}^2\to\Sp^{3\cdot2^{\mu-1}+1}$ of the projective plane $\mathbb{FP}^2$, where $\mathbb{F}=\R,\mathbb{C},\mathbb{H},\mathbb{O}$, respectively.

Conversely, for every conformal transformation $\phi$ of $\Sp^{3\cdot2^{\mu-1}+1}$, the composition $\phi\circ\psi_{\mathbb{F}}$ is a 2-Dupin isometric immersion.
\enthr
\bere
The remaining case $p=1$ of hypersurfaces in the above theorem corresponds to Pinkall's \cite{Pinkall1985} celebrated characterization of the well-known cyclides of Dupin as being always locally conformally congruent to rotational hypersurfaces with a sphere as profile.
\enre
The proof of Theorem \ref{omth} relies on Pinkall \cite{Pinkall1985ii} and Cecil-Jensen's \cite{Cecil1998} deep result that an irreducible 3-Dupin hypersurface must be Lie equivalent to a 3-isoparametric hypersurface, which in turn is a tube over the standard embedding $\psi_{\mathbb{F}}$ of the projective plane $\mathbb{FP}^2$, for $\mathbb{F}=\R,\mathbb{C},\mathbb{H},\mathbb{O}$, in $\Sp^4,\Sp^7,\Sp^{13},\Sp^{25}$, respectively, by Cartan's classification. Given a 2-Dupin submanifold $f:M^n\to\Sp^{n+p},\ p\geq2,$ our proof can be briefly outlined as follows.
\begin{enumerate}
\item We consider the Legendre lift $\lambda:UN_f^{n+p-1}\to\Lambda^{2\lp n+p\rp-1}$ of $f$, which will be a 3-Dupin submanifold in the sense of Lie sphere geometry.
\item We manage to prove that $\lambda$ is irreducible in the sense of Pinkall \cite{Pinkall1985}.
\item We then use Pinkall and Cecil-Jensen's aforementioned result coupled with Cartan's classification of 3-isoparametric hypersurfaces to conclude that $\lambda$ is Lie equivalent to a tube over some $\psi_{\mathbb{F}}$.
\item We use a characterization of Lie sphere transformations by Cecil-Chern \cite{Cecil1987} according to which any such transformation can be written as the composition of two M\"{o}bius transformations and some parallel transformation in the middle.
\item We finally show that if the above parallel transformation existed, then $f$ would be 2-unipotent, and therefore it would already be itself one of the standard embeddings $\psi_{\mathbb{F}}$, by Theorem \ref{dsth}.
\end{enumerate}

We now turn our attention to the broader class of 2-umbilical submanifolds of space forms. The fact mentioned above that the only principal curvatures which must be considered in checking the Dupin property are those of multiplicity one takes us one step away from the complete classification of these submanifolds too. In fact, as we will see in the next section, the two distinct principal curvatures of a 2-umbilical submanifold $f:M^n\to\Q^{n+p}_c$ of codimension $p>1$ have always the same multiplicity $m$, and in particular $M^n$ has even dimension $n=2m$. Therefore, every 2-umbilical submanifold that is not a surface nor a hypersurface is automatically 2-Dupin, and we are covered by Theorem \ref{omth}. In the case of a surface $f:M^2\to\Q^{2+p}_c$, there is not much to say, except that $p\leq2$, since the vector subspace of real symmetric $2\times2$ matrices is 3-dimensional and $f$ has no umbilical normal directions by definition. Thus, the only possibility remaining is that $f:M^n\to\Q^{n+1}_c$ is a hypersurface with a principal curvature $\kappa$ of multiplicity $n-1\geq2$. These are the well-known generalized $\lp n-1\rp$-cylinders for $\kappa=0$ and 1-parameter envelopes of hyperspheres (briefly, 1-PES) otherwise, which were studied by Dajczer-Florit-Tojeiro \cite{dajczer2005} and Asperti-Dajczer \cite{MR761994}, respectively. For the sake of completeness, we include a brief account of these results in Subsection \ref{1pes}, where we also offer an erratum to one of Asperti-Dajczer's results about the so-called special $k$-parameter envelopes of spheres (briefly, $k$-SPES).

Coupling Theorem \ref{omth} with the above considerations, we get the following complete classification of 2-umbilical submanifolds of space forms. Since the $k$-umbilical property is also invariant under conformal transformations between ambient spaces (because the traceless parts of the shape operators do not change under such a transformation), we present also this result just in the sphere.
\bethr\label{2uth}
Let $f:M^n\to\Sp^{n+p}$ be a 2-umbilical isometric immersion. Then, one of the following possibilities holds:
\begin{enumerate}[(i)]
\item $n=2^\mu$, $p=2^{\mu-1}+1$ for $\mu=2,3,4$ and there exists a conformal transformation $\phi$ of $\Sp^{3\cdot2^{\mu-1}+1}$ such that $\phi\circ f\lp M\rp$ is an open subset of the image of the standard embedding $\psi_{\mathbb{F}}:\mathbb{FP}^2\to\Sp^{3\cdot2^{\mu-1}+1}$ of the projective plane $\mathbb{FP}^2$, where $\mathbb{F}=\mathbb{C},\mathbb{H},\mathbb{O}$, respectively.
\item $n\geq3$, $p=1$ and $f$ is locally conformally congruent to a rotational hypersurface with a sphere as profile.
\item $n\geq3$, $p=1$ and, on each connected component of an open dense subset of $M^n$, $f$ is either a generalized $\lp n-1\rp$-cylinder or a 1-PES.
\item $n=2$, $p\leq2$ and $f$ is a surface free of umbilical normal directions.
\end{enumerate}
\enthr
Finally, we also study a weaker notion of 2-umbilicality that allows the existence of umbilical normal directions. An isometric immersion $f:M^n\to\tilde{M}^m$ is said to be \emph{weakly} $k$\emph{-umbilical} at $x\in M^n$ if it has at most $k$ distinct principal curvatures with respect to all $\lp x,\xi\rp\in UN_f^{m-1}$. A submanifold is called \emph{weakly} $k$\emph{-umbilical} if it is weakly $k$-umbilical at every point. In this case the principal curvatures $\kappa_1\leq\kappa_2\leq\dots\leq\kappa_n$ of $f$ with respect to $\lp x,\xi\rp$ are still continuous in the whole $UN_f^{m-1}$, but since they no longer have constant multiplicity, smoothness can be now ensured only on a open dense subset of $UN_f^{m-1}$. Note also that with this definition every weakly $\lp k-1\rp$-umbilical submanifold is also weakly $k$-umbilical. In particular, umbilical submanifolds are weakly 2-umbilical.
\bere
The above idea of weak $k$-umbilicality should not be confused with the concept of weak-umbilic points introduced by Moore \cite{moore1977}.
\enre
One might hope that this new concept gives rise to some other interesting types of submanifolds in space forms in the case $k=2$. However, we break this expectation by showing that the new possibilities that come out of this definition are not much distinct from the examples we have already gotten in Theorem \ref{2uth}. Before stating the result, first note that the standard embeddings in case (\emph{i}) of Theorem \ref{2uth} have flat normal bundle nowhere, while the rotational hypersurfaces, generalized cylinders and envelopes of hyperspheres in cases (\emph{ii}) and (\emph{iii}), being hypersurfaces, have flat normal bundle everywhere. It turns out that this duality is preserved in the weakly $2$-umbilical case. We show that allowing umbilical normal directions keeps cases (\emph{i}) and (\emph{ii}) essentially unaffected, as they remain true up to reduction of conformal codimension. Regarding case (\emph{iii}), the notion of a generalized cylinder actually makes sense in arbitrary codimension, while 1-PES are replaced by their natural generalization in higher codimension described in Subsection \ref{1pes}. Lastly, in case (\emph{iv}) all surfaces of arbitrary codimension are now clearly allowed. To avoid technicalities we shall assume $f$ in the statement to be umbilic free.
\bethr\label{w2uth}
Let $f:M^n\to\Sp^{n+p}$ be an umbilic-free weakly 2-umbilical isometric immersion. Then, one of the following possibilities holds:
\begin{enumerate}[(i)]
\item $R^\perp\neq0$ everywhere, $n=2^\mu$ for $\mu=2,3,4$ and there exist an umbilical sphere $\Sp^{3\cdot2^{\mu-1}+1}_c,\ c\geq1,$ of $\Sp^{2^\mu+p}$ such that $f\lp M\rp\subset\Sp^{3\cdot2^{\mu-1}+1}_c$ and a conformal transformation $\phi$ of $\Sp^{3\cdot2^{\mu-1}+1}_c$ such that $\theta_c^{-1}\circ\phi\circ f\lp M\rp$ is an open subset of the image of the standard embedding $\psi_{\mathbb{F}}:\mathbb{FP}^2\to\Sp^{3\cdot2^{\mu-1}+1}$ of the projective plane $\mathbb{FP}^2$, where $\mathbb{F}=\mathbb{C},\mathbb{H},\mathbb{O}$, respectively, and $\theta_c:\Sp^{3\cdot2^{\mu-1}+1}_c\to\Sp^{3\cdot2^{\mu-1}+1}$ is the similarity.
\item $R^\perp=0$ everywhere, $n\geq3$ and there exists an umbilical sphere $\Sp^{n+1}_c,\ c\geq1,$ of $\Sp^{n+p}$ such that $f\lp M\rp\subset\Sp^{n+1}_c$ and $f$ is locally conformally congruent to a rotational hypersurface in $\Sp^{n+1}_c$ with a sphere as profile.
\item $R^\perp=0$ everywhere, $n\geq3$ and, on each connected component of an open dense subset of $M^n$, $f$ either is a generalized $\lp n-1\rp$-cylinder or envelopes a 1-parameter congruence of hyperspheres.
\item $n=2$ and $f$ is an umbilic-free surface.
\end{enumerate}
\enthr
\noindent\small\emph{Acknowledgments.} The authors express deep thanks to E. Musso for the precious suggestion.

The second author is also looking forward to Rio 2016 and he hopes to see everybody there!
\normalsize
\section{Preliminaries}
Before proving our main theorems, we begin with some preliminary results on 2-umbilical submanifolds, Dupin submanifolds and 1-PES. To facilitate the exposition, we will divide this section into three subsections. The results in the first subsection derive from purely algebraic considerations on symmetric bilinear forms mirroring second fundamental forms of 2-umbilical submanifolds, while the second subsection consists of well-known results from Lie sphere geometry and its applications to the theory of Dupin submanifolds. The last subsection is devoted to the study of generalized $\lp n-1 \rp$-cylinders and  1-PES and their generalization to higher codimension.
\subsection{2-umbilical submanifolds} In this subsection, we derive two important results concerning 2-umbilical submanifolds. The first of these provides a simple (necessary and) sufficient condition for a 2-umbilical submanifold to be 2-unipotent and the second describes the spectral structure of the second fundamental form of 2-umbilical submanifolds.

In order to state the above lemmas we need to introduce some terminology. Let $V$ and $W$ be real vector spaces of finite dimension with positive definite inner products and let $\beta:V\times V\to W$ be a symmetric bilinear form. For any given $\xi\in W$, let $B_\xi:V\to V$ be the self-adjoint operator defined by
\beeqs
\la B_\xi X,Y\ra=\la\beta\lp X,Y\rp,\xi\ra.
\eneqs
Motivated by the situation of the second fundamental form of k-umbilical submanifolds, $\beta$ is said to be $k$\emph{-umbilical} if $B_\xi$ has exactly $k$ distinct eigenvalues $\kappa_i\lp\xi\rp,\ 1\leq i\leq k,$ for all $\xi\in W$ with $\left\|\xi\right\|=1$. It is called $k$\emph{-unipotent} if in addition either $k=1$ or $k\geq2$ and the eigenvalues $\kappa_i\lp\xi\rp=\kappa_i,\ 1\leq i\leq k,$ are constant.

In our target case of 2-umbilical symmetric bilinear forms, the next lemma shows that it suffices to check that all $B_\xi$'s share just one eigenvalue for $\beta$ to be 2-unipotent.
\beles
Let $\beta:V^n\times V^n\to W^p$ be a 2-umbilical symmetric bilinear form. If all $B_\xi,\ \left\|\xi\right\|=1,$ share a common eigenvalue $\kappa$, then $\beta$ is 2-unipotent.
\enles
\bepr
In the case $p=1$ there is nothing to prove. Otherwise, we first claim that $\kappa\neq0$. To see this, consider the vector $\mathcal{H}\in W^p$ defined by
\beeqs
\mathcal{H}=\frac{1}{n}\sum_{j=1}^n\beta\lp X_j,X_j\rp
\eneqs
in terms of an orthonormal basis $X_1,\dots,X_n$ of $V^n$.

The preceding expression implies that
\beeqs
n\la\mathcal{H},\xi\ra=\trac B_\xi
\eneqs
for any $\xi\in W^p$, hence
\beeqs
\trac B_\xi=0
\eneqs
for any unit vector $\xi\in W^p$ orthogonal to $\mathcal{H}$. Since $B_\xi$ has exactly two distinct eigenvalues, the above equation gives
\beeqs
\kappa\neq0,
\eneqs
as we wished. Then, since by assumption $\kappa$ must be a common eigenvalue of both $B_\xi$ and $B_{-\xi}=-B_\xi$, it follows that $\pm\kappa$ are their two eigenvalues.
\enpr
The above lemma has the following geometric counterpart.
\bele\label{upshl}
Let $f:M^n\to\Q^{n+p}_c$ be a 2-umbilical isometric immersion. If all shape operators $A_\xi,\ \lp x,\xi\rp\in UN_f^{n+p-1},$ share a common principal curvature $\kappa$, then $f$ is 2-unipotent.
\enle
\bepr
This is an immediate consequence of the previous lemma for $\beta=\alpha_x,\ x\in M^n.$
\enpr
Next, we show that 2-umbilicality imposes a very special spectral structure for symmetric bilinear forms. In particular, it will follow that, in codimension $p\geq2$, this phenomenon is unique to even dimensions.
\begin{lemma*}
Let $\beta:V^n\times V^n\to W^p,\ p\geq2,$ be a 2-umbilical symmetric bilinear form. Then, for every nonzero vector $\xi\in W^p$, we have
\beeq\label{eole}
\kappa_i\lp-\xi\rp=-\kappa_j\lp\xi\rp,\ 1\leq i\neq j\leq2,
\eneq
and thus
\beeq\label{eme}
m_1\lp\xi\rp=m_2\lp\xi\rp=:m.
\eneq
In particular, $n=2m$ is even.
\end{lemma*}
\bepr
Since $B_\xi$ has always two distinct eigenvalues, so does its traceless part $B_\xi^0$, whose eigenvalues are then $\pm\kappa_0\lp\xi\rp\neq0$. As $p\geq2$, we can walk from $\xi$ to $-\xi$ without stepping foot at the origin, so $\kappa_0$ never vanishes on the way and hence
\beeqs
\kappa_0\lp-\xi\rp=\kappa_0\lp\xi\rp,
\eneqs
which easily yields \eqref{eole}. Moreover, the fact that $E_{-\kappa_j}\lp-\xi\rp=E_{\kappa_j}\lp\xi\rp$, together with \eqref{eole}, clearly implies \eqref{eme}.
\enpr
Again we state the immediate geometric counterpart of the previous lemma.
\bele\label{edl}
Let $f:M^n\to\Q^{n+p}_c,\ p\geq2,$ be a 2-umbilical isometric immersion. Then, for every $\lp x,\xi\rp\in UN_f^{n+p-1}$, we have
\beeqs
\kappa_i\lp x,-\xi\rp=-\kappa_j\lp x,\xi\rp,\ 1\leq i\neq j\leq2,
\eneqs
and thus
\beeqs
m_1\lp x,\xi\rp=m_2\lp x,\xi\rp=:m.
\eneqs
In particular, $n=2m$ is even.
\enle
\bepr
Again a direct consequence of the preceding lemma for $\beta=\alpha_x,\ x\in M^n.$
\enpr
\bere\label{le2ex}
By relatively simple similar arguments the statements of Lemmas \ref{upshl} and \ref{edl} can be extended to weakly 2-umbilical submanifolds:
\item (i) The corresponding version of Lemma \ref{upshl} assumes that $f$ is umbilic-free and the shape operators $A_\xi$ share a common principal curvature $\kappa$ for all nonumbilical normal directions $\lp x,\xi\rp\in UN_f^{n+p-1}$. The conclusion is that either $\kappa=0$ and the first normal space $N^f_1\lp x\rp=\Span\lb\mathcal{H}\lp x\rp\rb$ has dimension 1 for all $x\in M^n$, in particular, the second fundamental form of $f$ is diagonalizable and $f$ has flat normal bundle everywhere by the Ricci equation, or $\kappa\neq0$ and $f$ is 2-unipotent;
\item (ii) The assertions in Lemma \ref{edl} remain true for every nonumbilical normal direction $\lp x,\xi\rp\in UN_f^{n+p-1}$ at every point $x\in M^n$ such that $R^\perp\lp x\rp\neq0$ in the weakly 2-umbilical case, in particular, $n=2m$ is even if $f$ has nonflat normal bundle somewhere.
\enre
\subsection{Lie sphere geometry and Dupin submanifolds}\label{lsgds}
Lie \cite{Lie1872} introduced the geometry of oriented spheres in his dissertation, published as a paper in \emph{Mathematische Annalen} in 1872. Sphere geometry was also prominent in his study of contact transformations (Lie-Scheffers \cite{MR0460049}) and in Volume III of Blaschke's book on differential geometry published in 1929.

Lie sphere geometry has become a valuable tool in the study of Dupin submanifolds, beginning with Pinkall's \cite{Pinkall1981} dissertation in 1981. In this subsection, we will present a series of results that will illustrate this value. The first result is an important characterization of the Lie sphere group in terms of M\"{o}bius and parallel transformations due to Cecil and Chern \cite{Cecil1987}. The next couple of results provide the basic link between Dupin submanifolds as defined in the introduction of this article and those in the context of Legendre submanifolds by means of the so-called Legendre lifts. This connection will allow us to deduce several properties of our Dupin submanifolds out of the ones satisfied by their Legendre lifts. It will follow, for example, that the only principal curvatures which must be considered in checking the Dupin property are those of multiplicity one, and also that our class of Dupin submanifolds is preserved by conformal transformations of and between ambient space forms.

Finally, we will discuss the classification of proper Dupin submanifolds with three curvature spheres by Cecil and Jensen \cite{Cecil1998}. A crucial assumption in their result is the irreducibility of the Dupin submanifold in the sense of Pinkall \cite{Pinkall1985}. In the next section, we will manage to prove that the Legendre lifts arising from our Dupin submanifolds always satisfy such condition, and for this purpose stating Pinkall's \cite{Pinkall1981} simple Lie sphere geometric criterion for reducibility will come in handy.

Before and during our presentation of the above results, we will be briefly reviewing some important concepts from Lie sphere geometry as we need them. A detailed account on the subject can be found in \cite{MR2361414}.

Start with the real vector space $\R^{d+3}_2$ endowed with the inner product of signature $\lp-1,1,\dots,1,-1\rp$ and define an equivalence relation on $\R^{d+3}_2-\lb0\rb$ by setting $x\simeq y$ if $x=ty$ for some nonzero real number $t$. We denote the equivalence class determined by a vector $x$ by $\lrb x\rrb$. Projective space $\mathbb{P}^{d+2}$ is the set of such equivalence classes, and it can naturally be identified with the space of all lines through the origin in $\R^{d+3}_2$. Let $Q^{d+1}$ be the quadric in $\mathbb{P}^{d+2}$ given in homogeneous coordinates by the equation
\beeq\label{lqeq}
\la x,x\ra=-x_1^2+x_2^2+\cdots+x_{d+2}^2-x_{d+3}^2=0.
\eneq
The manifold $Q^{d+1}$ is called the \emph{Lie quadric},  and the scalar product determined by the quadratic form in \eqref{lqeq} is called the \emph{Lie metric} or \emph{Lie scalar product}. We will let $\lb e_1,\dots,e_{d+3}\rb$ denote the standard orthonormal basis for $\R^{d+3}_2$. It turns out that points on $Q^{d+1}$ correspond to the set of oriented hyperspheres and point spheres (hyperspheres with radius 0) in $\Sp^d_c$, or in general to the set of oriented extrinsic hyperspheres and point spheres in $\Q^d_c\cup\lb\infty\rb$ for the Euclidean or a hyperbolic space form $\Q^d_c,\ c\leq0,$ where the \emph{improper point} $\infty$ is identified with the pole $\sigma$ of the stereographic projection $\pi_c$ given in equation \eqref{sgpr}. Actually, from the point of view of Klein's Erlangen Program, all three spherical, Euclidean and hyperbolic geometries are subgeometries of Lie sphere geometry (see \cite{MR2361414} for the details). On the other hand, in some ways it is simpler to use the sphere $\Sp^d$ rather than $\R^d$ or $\Hy^d$ as the base space for the study of M\"{o}bius or Lie sphere geometry, since this avoids the use of stereographic projection and the need to refer to an improper point or to distinguish between spheres, planes and horospheres. Indeed, the above correspondence between points in $Q^{d+1}$ and oriented spheres in the base space is much simpler in $\Sp^d$ than in $\R^d$ and $\Hy^d$. It is given simply by
\beeq\label{lscr}
S\lp x,r\rp\leftrightarrow\lrb\lp\cos r,x,\sin r\rp\rrb,
\eneq
where $S\lp x,r\rp$ denotes the sphere in $\Sp^d$ with center $x$ and signed radius $r$. In particular, the point sphere $x=S\lp x,0\rp$ in $\Sp^d$ corresponds to the point $\lrb\lp1,x,0\rp\rrb$ in $Q^{d+1}$. It is an important fact that the Lie quadric contains projective lines but no linear subspaces of higher dimension.

The set of oriented spheres in $\Sp^d$ corresponding to the points on a line on $Q^{d+1}$ forms a so-called \emph{parabolic pencil} of spheres. Each parabolic pencil contains exactly one point sphere $k_1$ and one great sphere $k_2$. The pencil consists of all oriented hyperspheres in oriented contact with $k_1$ and $k_2$.

A \emph{Lie sphere transformation} is a projective transformation of $\mathbb{P}^{d+2}$ which takes $Q^{d+1}$ to itself. In terms of the geometry of $\Sp^d$, a Lie sphere transformation maps Lie spheres to Lie spheres. (Here the term ``Lie sphere'' includes oriented spheres and point spheres.) Furthermore, since it is projective, a Lie sphere transformation maps lines on $Q^{d+1}$ to lines on $Q^{d+1}$. Thus, it preserves oriented contact of spheres in $\Sp^d$. The group $G$ of Lie sphere transformations is isomorphic to $O\lp d+1,2\rp/\lb\pm I\rb$, where $O\lp d+1,2\rp$ is the group of orthogonal transformations of $\R^{d+3}_2$. Pinkall's \cite{Pinkall1984}-\cite{Pinkall1984ii} so-called ``fundamental theorem of Lie sphere geometry'' states that any line preserving diffeomorphism of $Q^{d+1}$ is the restriction to $Q^{d+1}$ of a projective transformation, that is, a transformation of the space of oriented spheres which preserves oriented contact must be a Lie sphere transformation.

Recall that a linear transformation $L\in GL\lp d+1\rp$ induces a projective transformation $P\lp L\rp$ on $\mathbb{P}^d$ defined by $P\lp L\rp\lrb x\rrb=\lrb Lx\rrb$. The map $P$ is a homomorphism of $GL\lp d+1\rp$ onto the group $PGL\lp d\rp$ of projective transformations of $\mathbb{P}^d$. It is well known (see, for example, Samuel \cite{MR960691}) that the kernel of $P$ is the group of all nonzero scalar multiples of the identity transformation $I$. A M\"{o}bius transformation is a transformation on the space of unoriented spheres, i.e., the space of projective classes of spacelike vectors in $\R^{d+2}_1$. Such a transformation also takes lightlike vectors to lightlike vectors, and so it induces a conformal diffeomorphism of the sphere $\Sp^d$ (regarded as the M\"{o}bius quadric in $\mathbb{P}^{d+1}$) onto itself.

It is well known that the group of conformal diffeomorphisms of the sphere is precisely the M\"{o}bius group. Now, each M\"{o}bius transformation naturally induces two Lie sphere transformations on the space $Q^{d+1}$ of oriented spheres. Specifically, if $L$ is in $O\lp d+1,1\rp$, then we can extend $L$ to a transformation $\bar{L}$ in $O\lp d+1,2\rp$ by setting $\bar{L}=L$ on $\R^{d+2}_1$ and $\bar{L}\lp e_{d+3}\rp=\pm e_{d+3}$, where the sign determines the orientation of every oriented sphere. Thus, the group of Lie transformations induced from M\"{o}bius transformations consists of those Lie transformations that map $\lrb e_{d+3}\rrb$ to itself. Since such a transformation must also take $e_{d+3}^\perp$ to itself, this is precisely the group of Lie transformations which take point spheres to point spheres.

When working in the context of Lie sphere geometry, it is a common abuse of language to refer to these transformations also as ``M\"{o}bius transformations''. Besides M\"{o}bius transformations, an important family of Lie transformations is given by the so-called parallel transformations $P_t$, which add $t$ to the signed radius of each sphere while keeping the center fixed. For each of the spherical, Euclidean and hyperbolic geometries, there is a corresponding one-parameter family of parallel transformations. \emph{Spherical parallel transformation} $P_t$ is accomplished by the following transformation in $O\lp d+1,2\rp$:
\beeq\label{spteq}
\begin{aligned}
P_te_1 & =\cos t\,e_1+\sin t\,e_{d+3},\\
P_te_{d+3} & =-\sin t\,e_1+\cos t\,e_{d+3},\\
P_te_i & =e_i,\quad 2\leq i\leq d+2.
\end{aligned}
\eneq
On the other hand, \emph{Euclidean parallel transformation} $P_t$ is accomplished in terms of matrix representation with respect to the standard orthonormal basis by
\beeq
P_t=
  \begin{bmatrix}
    1-\lp t^2/2\rp & -t^2/2 & 0 & -t \\
    t^2/2 & 1+\lp t^2/2\rp & 0 & t \\
    0 & 0 & I & 0 \\
    t & t & 0 & 1
  \end{bmatrix}.
\eneq
Lastly, \emph{hyperbolic parallel transformation} is accomplished by the following transformation:
\beeq
\begin{aligned}
P_te_i & =e_i,\quad i=1,3,\dots,d+2,\\
P_te_2 & =\cosh t\,e_2+\sinh t\,e_{d+3},\\
P_te_{d+3} & =\sinh t\,e_2+\cosh t\,e_{d+3}.
\end{aligned}
\eneq
The following theorem demonstrates the important role played by parallel transformations in generating the Lie sphere group (see Cecil-Chern \cite{Cecil1987}).
\bethr[Cecil-Chern \cite{Cecil1987}]\label{thcc}
Any Lie sphere transformation $\gamma$ can be written as
\beeqs
\gamma=\Phi_1P_t\Phi_2,
\eneqs
where $\Phi_1$ and $\Phi_2$ are M\"{o}bius transformations and $P_t$ is some Euclidean, spherical or hyperbolic parallel transformation.
\enthr

Let us now develop the framework necessary to study submanifolds within the context of Lie sphere geometry. The manifold $\Lambda^{2d-1}$ of projective lines on the Lie quadric $Q^{d+1}$ can be identified with the tangent bundle $T_1\Sp^d$ of $\Sp^d$ and thus borrows from the latter a contact structure, i.e., a globally defined 1-form $\omega$ such that $\omega\wedge\lp d\omega\rp^{d-1}\neq0$ on $\Lambda^{2d-1}$. This gives rise to a codimension one distribution $D$ on $\Lambda^{2d-1}$ that has integral submanifolds of dimension $d-1$, but none of higher dimension. These integral submanifolds are called Legendre submanifolds. Any submanifold of a real space-form $\Sp^d,\ \R^d$ or $\Hy^d$ naturally induces a Legendre submanifold, and thus Lie sphere geometry can be used to analyze submanifolds in these spaces. This has been particularly effective in the classification of Dupin submanifolds.

To see how a submanifold in a space form gives rise to a Legendre submanifold, first suppose that $f:M^n\to\Sp^{n+1}$ is an immersed oriented hypersurface with field of unit normals $\xi:M^n\to\Sp^{n+1}$. The induced Legendre submanifold is then given by the map $\lambda:M^n\to\Lambda^{2n+1}$ defined by
\beeqs
\lambda\lp x\rp=\lrb Y_1\lp x\rp,Y_{n+4}\lp x\rp\rrb,
\eneqs
where
\beeqs
Y_1\lp x\rp=\lp 1,f\lp x\rp,0\rp,\quad Y_{n+4}\lp x\rp=\lp0,\xi\lp x\rp,1\rp.
\eneqs
The map $\lambda$ is called the \emph{Legendre lift} of the immersion $f$ with field of unit normals $\xi$. Next, we handle the case of a submanifold $f:M^n\to\Sp^{n+p}$ of codimension $p$ greater than one. The \emph{Legendre lift} of $f$ is the map $\lambda:UN_f^{n+p-1}\to\Lambda^{2\lp n+p\rp-1}$ defined by
\beeqs
\lambda\lp x,\xi\rp=\lrb Y_1\lp x\rp,Y_{n+p+3}\lp x,\xi\rp\rrb,
\eneqs
where
\beeq\label{lldf}
Y_1\lp x\rp=\lp1,f\lp x\rp,0\rp,\quad Y_{n+p+3}\lp x,\xi\rp=\lp0,\xi,1\rp.
\eneq
Geometrically, $\lambda\lp x,\xi\rp$ is the line on the quadric $Q^{n+p+1}$ corresponding to the parabolic pencil of spheres in $\Sp^{n+p}$ in oriented contact at the contact element $\lp f\lp x\rp,\xi\rp\in T_1\Sp^{n+p}$. The situation for submanifolds of $\R^{n+p}$ or $\Hy^{n+p}$ is similar and we shall not treat it here (cf. Cecil \cite{MR2361414}). Now suppose that $\lambda:N^{d-1}\to\Lambda^{2d-1}$ is an arbitrary Legendre submanifold. Then, it is always possible to parametrize $\lambda$ by the point sphere map $\lrb Y_1\rrb$ and the great sphere map $\lrb Y_{d+3}\rrb$ given by
\beeq\label{psgsmeq}
Y_1=\lp1,f,0\rp,\quad Y_{d+3}=\lp0,\xi,1\rp.
\eneq
This defines two maps $f$ and $\xi$ from $N^{d-1}$ to $\Sp^d$, which we call the \emph{spherical projection} and \emph{spherical field of unit normals}, respectively. Both $f$ and $\xi$ are smooth maps, but neither need to be an immersion or even have constant rank. The Legendre submanifold induced by an oriented hypersurface in $\Sp^d$ is the special case where the spherical projection $f$ is an immersion, i.e., $f$ has constant rank $d-1$ on $N^{d-1}$. In the case of the Legendre submanifold induced by a submanifold $f:M^n\to\Sp^{n+p}$, the spherical projection $f:UN_f^{n+p-1}\to\Sp^{n+p}$ given by $f\lp x,\xi\rp=f\lp x\rp$ is a submersion onto the image of $f$, i.e., it has constant rank $n$. The notions of projection and field of unit normals also have their Euclidean and hyperbolic versions.
\bere
The jump from 1 to $d+3,\ d=n+p,$ in the indices of the two $Y_i$'s defining the Legendre lift may seem a bit messy, but the notation is traditional from the method of moving frames in the context of Lie sphere geometry, where $Y_1$ and $Y_{d+3}$ are just the first and last vectors of an ordered set called a Lie frame (refer to Cecil-Chern \cite{Cecil1987}).
\enre
We now discuss the definition of a curvature sphere of an isometric immersion $f:M^n\to\Q^{n+p}_c$, which is closely related to focal points of $f$, that is, singular points of the \emph{tube} $f_t:UN_f^{n+p-1}\to\Q^{n+p}_c$ \emph{of radius} $t$ \emph{around} $f$ given by
\beeqs
\begin{aligned}
f_t\lp x,\xi\rp & =\cos\lp\sqrt{c}t\rp f\lp x\rp+\frac{\sin\lp\sqrt{c}t\rp}{\sqrt{c}}\xi,\quad\text{if }c>0,\\
f_t\lp x,\xi\rp & =f\lp x\rp+t\xi,\quad\text{ if }c=0,\\
f_t\lp x,\xi\rp & =\cosh\lp\sqrt{-c}t\rp f\lp x\rp+\frac{\sinh\lp\sqrt{-c}t\rp}{\sqrt{-c}}\xi,\quad\text{if }c<0.
\end{aligned}
\eneqs
Note that if $f:M^n\to\Q^{n+1}_c$ is an oriented hypersurface with unit normal vector field $\xi$, then $UN_f^n=\lb\lp x,\xi\lp x\rp\rp:x\in M^n\rb$ and the above tube is essentially the \emph{parallel hypersurface} $f_t:M^n\to\Q^{n+1}_c$ \emph{of} $f$ \emph{at distance} $t$ defined by
\beeqs
\begin{aligned}
f_t\lp x\rp & =\cos\lp\sqrt{c}t\rp f\lp x\rp+\frac{\sin\lp\sqrt{c}t\rp}{\sqrt{c}}\xi\lp x\rp,\quad\text{if }c>0,\\
f_t\lp x\rp & =f\lp x\rp+t\xi\lp x\rp,\quad\text{ if }c=0,\\
f_t\lp x\rp & =\cosh\lp\sqrt{-c}t\rp f\lp x\rp+\frac{\sinh\lp\sqrt{-c}t\rp}{\sqrt{-c}}\xi\lp x\rp,\quad\text{if }c<0.
\end{aligned}
\eneqs
For simplicity we assume $c=1$, so
\beeq\label{tbeq}
f_t\lp x,\xi\rp=\cos t\,f\lp x\rp+\sin t\,\xi.
\eneq
Geometrically, one thinks of focal points as points where nearby normal geodesics intersect. It is well known that the location of focal points is related to the principal curvatures. Specifically, if $\lp X,-f_*A_\xi X+\eta\rp\in T_{\lp x,\xi\rp}UN_f$, then we have
\beeqs
\begin{split}
{f_t}_*\lp X,-f_*A_\xi X+\eta\rp=\cos t\,f_*X+\sin t\,\lp-f_*A_\xi X+\eta\rp\\
=f_*\lp\cos t\,X-\sin t\,A_\xi X\rp+\sin t\,\eta.
\end{split}
\eneqs
Thus, ${f_t}_*\lp X,-f_*A_\xi X+\eta\rp=0$ for $\lp X,-f_*A_\xi X+\eta\rp\neq0$ if and only if either $p\geq2,\ t=k\pi,\ k\in\mathbb{Z},$ and $X=0$ or $t\neq k\pi,\ k\in\mathbb{Z},\ \eta=0,\ \cot t$ is a principal curvature of $f$ at $\lp x,\xi\rp$ and $X$ is a corresponding principal vector. Hence, $f_{k\pi}\lp x,\xi\rp=\lp-1\rp^kf\lp x\rp,\ k\in\mathbb{Z},$ is a focal point of $f$ at $\lp x,\xi\rp$ of multiplicity $p-1$ if $p\geq2$ and $y=f_t\lp x,\xi\rp$ is a focal point of $f$ at $\lp x,\xi\rp$ of multiplicity $m$ if and only if $\cot t$ is a principal curvature of multiplicity $m$ at $\lp x,\xi\rp$. If $p\geq2$, it is convenient to introduce an extra principal curvature $\kappa=\cot0=\infty$ of multiplicity $p-1$ associated to the focal points $\pm f\lp x\rp$. Then each principal curvature
\beeqs
\kappa=\cot t,\quad 0\leq t<\pi,
\eneqs
produces two distinct antipodal focal points on the normal geodesic with parameter values $t$ and $t+\pi$. The oriented hypersphere centered at a focal point $y=f_t\lp x,\xi\rp$ and in oriented contact with $f\lp M\rp$ at $\lp f\lp x\rp,\xi\rp$ is called a \emph{curvature sphere} of $f$ at $\lp x,\xi\rp$. The two antipodal focal points determined by $\kappa$ are the two centers of the corresponding curvature sphere. Thus, the correspondence between principal curvatures and curvature spheres is bijective.

The \emph{multiplicity} of the curvature sphere is by definition equal to the multiplicity of the corresponding principal curvature. The \emph{principal space} corresponding to the curvature sphere is by definition equal to the horizontal lift of the principal space of the corresponding principal curvature $\kappa=\cot t$ for $0<t<\pi$ or $\lb\lp0,\eta\rp:\eta\in N_fM\lp x\rp,\eta\cdot\xi=0\rb$ for $t=0$. We now consider these ideas as they apply to the Legendre lift of the isometric immersion $f:M^n\to\Sp^{n+p}$. As in equation \eqref{lldf}, we have $\lambda=\lrb Y_1,Y_{n+p+3}\rrb$, where
\beeqs
Y_1=\lp1,f,0\rp,\quad Y_{n+p+3}=\lp0,\xi,1\rp.
\eneqs
For each $\lp x,\xi\rp\in UN_f^{n+p-1}$, the points on the line $\lambda\lp x,\xi\rp$ can be parametrized as
\beeqs
\lrb K_t\lp x,\xi\rp\rrb=\lrb\cos t\,Y_1\lp x\rp+\sin t\,Y_{n+p+3}\lp x,\xi\rp\rrb=\lrb\lp\cos t,f_t\lp x,\xi\rp,\sin t\rp\rrb,
\eneqs
where $f_t$ is given in equation \eqref{tbeq}. By equation \eqref{lscr}, the point $\lrb K_t\lp x,\xi\rp\rrb$ in $Q^{n+p+1}$ corresponds to the oriented sphere in $\Sp^{n+p}$ with center $f_t\lp x,\xi\rp$ and signed radius $t$. This sphere is in oriented contact with the submanifold $f\lp M\rp$ at $\lp f\lp x\rp,\xi\rp$. Given a tangent vector $\lp X,-f_*A_\xi X+\eta\rp\in T_{\lp x,\xi\rp}UN_f$, we have
\beeqs
{K_t}_*\lp X,-f_*A_\xi X+\eta\rp=\lp0,{f_t}_*\lp X,-f_*A_\xi X+\eta\rp,0\rp.
\eneqs
Thus, ${K_t}_*\lp X,-f_*A_\xi X+\eta\rp=\lp0,0,0\rp$ if and only if ${f_t}_*\lp X,-f_*A_\xi X+\eta\rp=0$, i.e., $y=f_t\lp x,\xi\rp$ is a focal point of $f$ at $\lp x,\xi\rp$. Hence, we have that the point $\lrb K_t\lp x,\xi\rp\rrb$ in $Q^{n+p+1}$ corresponds to a curvature sphere of the isometric immersion $f$ at $\lp x,\xi\rp$ if and only if
\beeq\label{cslm}
{K_t}_*\lp X,-f_*A_\xi X+\eta\rp=\lp0,0,0\rp
\eneq
for some nonzero vector $\lp X,-f_*A_\xi X+\eta\rp\in T_{\lp x,\xi\rp}UN_f$. This viewpoint allows us to extend the definition of curvature spheres and $k$-umbilicality to an arbitrary Legendre submanifold in the obvious way.
\bere
Note that if we use a different parametrization
\beeqs
\lrb\tilde{K}_t\lp x,\xi\rp\rrb=\lrb\rho\lp x,\xi\rp K_t\lp x,\xi\rp\rrb
\eneqs
of $\lrb K_t\lp x,\xi\rp\rrb$, then condition \eqref{cslm} for $\tilde{K_t}$ is that $\tilde{K_t}_*\lp X,-f_*A_\xi X+\eta\rp$ is collinear with $\tilde{K_t}$.
\enre
One can show that the notion of curvature sphere is invariant under Lie sphere transformations. Let $\lambda:N^{d-1}\to\Lambda^{2d-1}$ be a Legendre submanifold parametrized by $\lambda=\lrb Z_1,Z_{d+3}\rrb$. Suppose $\gamma=P\lp L\rp$ is the Lie sphere transformation induced by an orthogonal transformation $L$ in the group $O\lp d+1,2\rp$. Since $L$ is orthogonal, it is easy to check that the map defined by $\lrb LZ_1,LZ_{d+3}\rrb$ is also a Legendre submanifold, which we will denote by $\gamma\lambda:N^{d-1}\to\Lambda^{2d-1}$. The Legendre submanifolds $\lambda$ and $\gamma\lambda$ are said to be \emph{Lie equivalent}. In terms of space form geometry, suppose that $f:M^n\to\Q^{n+p}_c$ and $g:M^n\to\Q^{n+p}_c$ are two isometric immersions. We say that $f$ and $g$ are \emph{Lie equivalent} if their Legendre lifts are Lie equivalent. It is easy to see that, for $\lambda$ and $\gamma$ as above, the point $\lrb K\rrb$ on the line $\lambda\lp x\rp$ is a curvature sphere of $\lambda$ at $x$ if and only if the point $\gamma\lrb K\rrb$ is a curvature sphere of the Legendre submanifold $\gamma\lambda$ at $x$. Furthermore, the principal spaces corresponding to $\lrb K\rrb$ and $\gamma\lrb K\rrb$ are identical. By Theorem \ref{thcc}, two important special cases are when the Lie sphere transformation $\gamma$ is a M\"{o}bius or a parallel transformation. In the first case, suppose that $\gamma=\Phi$ is induced by a conformal diffeomorphism $\phi:\bar{\Q}^{n+p}_c\to\bar{\Q}^{n+p}_{\tilde{c}}$, where
\beeqs
\bar{\Q}^{n+p}_c=
\begin{cases}
\Sp^{n+p}_c-\lb\sigma\rb,\quad\text{if }c>0,\\
\Q^{n+p}_c,\quad\text{if }c\leq0.
\end{cases}
\eneqs
If $f:M^n\to\bar{\Q}^{n+p}_c$ is an isometric immersion and $\tilde{f}=\phi\circ f:M^n\to\bar{\Q}^{n+p}_{\tilde{c}}$, then the Legendre lifts $\lambda:UN_f^{n+p-1}\to\Lambda^{2\lp n+p\rp-1}$ of $f$ and $\tilde{\lambda}:UN_{\tilde{f}}^{n+p-1}\to\Lambda^{2\lp n+p\rp-1}$ of $\tilde{f}$ are \emph{M\"{o}bius conjugated} in the sense that
\beeqs
\tilde{\lambda}=\Phi\lambda\bar{\phi}^{-1},
\eneqs
where $\bar{\phi}:UN_f^{n+p-1}\to UN_{\tilde{f}}^{n+p-1}$ is the unit normal bundle isometry induced by $\phi$. In the second case, suppose for simplicity that $\gamma=P_t$ is the spherical parallel transformation given in equation \eqref{spteq}. Recall that $P_t$ has the effect of adding $t$ to the signed radius of each sphere in $\Sp^d$ while keeping the center fixed. Suppose that $\lambda:N^{d-1}\to\Lambda^{2d-1}$ is a Legendre submanifold parametrized by the point sphere and great sphere maps $\lb Y_1,Y_{d+3}\rb$, as in equation \eqref{psgsmeq}. Then $P_t\lambda=\lrb W_1,W_{d+3}\rrb$, where
\beeqs
W_1=P_tY_1=\lp\cos t,f,\sin t\rp,\quad W_{d+3}=P_tY_{d+3}=\lp-\sin t,\xi,\cos t\rp.
\eneqs
Note that $W_1$ and $W_{d+3}$ are not the point sphere and great sphere maps for $P_t\lambda$. Solving for the point sphere map $Z_1$ and the great sphere map $Z_{d+3}$ of $P_t\lambda$, we find
\beeqs
\begin{aligned}
Z_1 & =\cos t\,W_1-\sin t\,W_{d+3}=\lp1,\cos t\,f-\sin t\,\xi,0\rp,\\
Z_{d+3} & =\sin t\,W_1+\cos t\,W_{d+3}=\lp0,\sin t\,f+\cos t\,\xi,1\rp.
\end{aligned}
\eneqs
From this, we see that $P_t\lambda$ has spherical projection and spherical unit normal field given, respectively, by
\beeq\label{spsunfpeq}
\begin{aligned}
f_{-t} & =\cos t\,f-\sin t\,\xi=\cos\lp-t\rp f+\sin\lp-t\rp\xi,\\
\xi_{-t} & =\sin t\,f+\cos t\,\xi=-\sin\lp-t\rp f+\cos\lp-t\rp\xi.
\end{aligned}
\eneq
The minus sign occurs because $P_t$ takes a sphere with center $f_{-t}\lp x\rp$ and radius $-t$ to the point sphere $f_{-t}\lp x\rp$. We call $P_t\lambda$ a \emph{parallel submanifold} of $\lambda$. Formula \eqref{spsunfpeq} shows the close correspondence between these parallel submanifolds and the parallel hypersurfaces to or the tubes $f_t$ around $f$, in the case where $f:M^n\to\Sp^{n+p}$ is an isometric immersion with $p=1$ or $p\geq2$, respectively. The spherical projection $f_t$ has singularities at the focal points of $f$, but the parallel submanifold $P_t\lambda$ is still a smooth submanifold of $\Lambda^{2d-1}$.

An important theorem, due to Pinkall \cite{Pinkall1985}, shows that the number of these singularities is bounded for each $x\in M^{d-1}$, which allows us to obtain several consequences by passing to a parallel submanifold, if necessary, and then applying well-known results concerning immersed hypersurfaces in $\Sp^d$. Among them are the important facts that if the dimension $m$ of the principal space $T$ of a given curvature sphere $K$ is constant on an open subset $U$ of $M^{d-1}$, then the principal distribution $T$ is integrable on $U$, and if $m>1$, then the curvature sphere map $K$ is constant along the leaves of the principal foliation $T$.

A connected submanifold $S$ of $M^{d-1}$ is called a \emph{curvature submanifold} if at each $x\in S$, the tangent space $T_xS$ is equal to some principal space $T$. For example, if $\dim T$ is constant on an open subset $U$ of $M^{d-1}$, then each leaf of the principal foliation $T$ is a curvature submanifold on $U$. Curvature submanifolds are plentiful, since the results of Reckziegel \cite{Reckziegel1979} and Singley \cite{singley1975} imply that there is an open dense subset $\Omega$ of $M^{d-1}$ on which the multiplicities of the curvature spheres are locally constant. On $\Omega$, each leaf of each principal foliation is a curvature submanifold. On the other hand, it is also possible to have a curvature submanifold $S$ which is not a leaf of a principal foliation, because the multiplicity of the corresponding curvature sphere is not constant on a neighborhood of $S$. The fact that curvature spheres of multiplicity $m\geq2$ are constant along the leaves of their corresponding principal foliations can be generalized to the fact that curvature spheres are always constant along corresponding curvature submanifolds of dimension $m\geq2$. The proof is obtained by invoking a theorem of Ryan \cite{ryan1969}.

The notion of Dupin submanifold is generalized to the context of Lie sphere geometry by calling a Legendre submanifold \emph{Dupin} if along each curvature submanifold, the corresponding curvature sphere is constant \cite{MR2361414}. A Dupin submanifold,
\beeqs
\lambda:M^{d-1}\to\Lambda^{2d-1},
\eneqs
is said to be \emph{proper Dupin} if the number $k$ of distinct curvature spheres is constant on $M^{d-1}$. In this case, we also refer to $\lambda$ as a $k$-Dupin submanifold. Of course, Legendre lifts of Dupin hypersurfaces in $\Sp^{n+1}$ are Dupin in the sense defined here. Indeed, we have the following basic link between the two definitions.
\bethr[Cecil \cite{MR2361414}]\label{cecil1}
Let $f:M^n\to\Sp^{n+1}$ be a hypersurface and $\lambda:M^n\to\Lambda^{2n+1}$ its Legendre lift. Let $Y_1$ and $Y_{n+4}$ be the point sphere and great sphere maps of $\lambda$, respectively. Then the curvature spheres of $\lambda$ at a point $x\in M^n$ are
\beeqs
\lrb K_i\rrb=\lrb\kappa_iY_1+Y_{n+4}\rrb,\ 1\leq i\leq k,
\eneqs
where $\kappa_1,\dots,\kappa_k$ are the distinct principal curvatures at $x$ of $f$. The principal space of $\lambda$ corresponding to the curvature sphere $\lrb K_i\rrb$ equals the principal space of $f$ corresponding to the principal curvature $\kappa_i$. In particular, the multiplicity of $\lrb K_i\rrb$ equals the multiplicity of $\kappa_i$.
Moreover, $\lrb K_i\rrb$ is Dupin if and only if $\kappa_i$ is Dupin.
\enthr
\bepr
This follows immediately from the fact that the point
\beeqs
\lrb K\rrb=\lrb\kappa Y_1+Y_{n+4}\rrb
\eneqs
is a curvature sphere of $\lambda$ at $x$ with corresponding principal vector $X\in T_xM$ if and only if
\beeqs
\begin{split}
\kappa{Y_1}_*X+{Y_{n+4}}_*X=\kappa\lp0,f_*X,0\rp+\lp0,-f_*AX,0\rp\\
=\lp0,-f_*\lp AX-\kappa X\rp,0\rp=0,
\end{split}
\eneqs
in which case
\beeqs
K_*X=\kappa_*X\,Y_1+\kappa{Y_1}_*X+{Y_{n+4}}_*X=\kappa_*XY_1
\eneqs
is collinear with $K$ if and only if $\kappa_*X=0$, as desired.
\enpr
Pinkall's \cite{Pinkall1986} definition of Dupin submanifolds in higher codimension can also be expressed in terms of their Legendre lifts.
\bethr[Cecil \cite{MR2361414}]\label{cecil2}
Let $f:M^n\to\Sp^{n+p},\ p\geq2,$ be an isometric immersion and $\lambda:UN_f^{n+p-1}\to\Lambda^{2\lp n+p\rp-1}$ its Legendre lift. Let $Y_1$ and $Y_{n+p+3}$ be the point sphere and great sphere maps of $\lambda$, respectively. Then the curvature spheres of $\lambda$ at a point $\lp x,\xi\rp\in UN_f^{n+p-1}$ are
\beeq\label{cspc}
\lrb K_i\rrb=\lrb\kappa_iY_1+Y_{n+p+3}\rrb,\ 1\leq i\leq k+1,
\eneq
where $\kappa_1,\dots,\kappa_k$ are the distinct principal curvatures at $\lp x,\xi\rp$ of $f$ and $\kappa_{k+1}=\infty$. For $1\leq i\leq k$, the principal space of $\lambda$ corresponding to the curvature sphere $\lrb K_i\rrb$ is given by
\beeqs
T_i=\lb\lp X,-f_*A_\xi X\rp:X\in E_i\rb,
\eneqs
where $E_i$ is the principal space of $f$ corresponding to the principal curvature $\kappa_i$, while the principal space corresponding to $\lrb K_{k+1}\rrb$ is
\beeq\label{psps}
T_{k+1}\lp x,\xi\rp=\lb\lp0,\eta\rp:\eta\in N_fM\lp x\rp,\eta\cdot\xi=0\rb.
\eneq
In particular, the multiplicity of $\lrb K_i\rrb$ equals the multiplicity of $\kappa_i$ for $1\leq i\leq k$, while the multiplicity of $\lrb K_{k+1}\rrb$ is $p-1$. Moreover, for $1\leq i\leq k,\ \lrb K_i\rrb$ is Dupin if and only if $\kappa_i$ is Dupin, while $\lrb K_{k+1}\rrb$ is always Dupin.
\enthr
\bepr
Note that
\beeqs
\lrb K_{k+1}\rrb=\lrb\infty Y_1+Y_{n+p+3}\rrb=\lrb Y_1\rrb
\eneqs
is simply the point sphere map of $\lambda$, which is clearly a curvature sphere with corresponding principal space $T_{k+1}$ given by \eqref{psps} and multiplicity $p-1$, since $p\geq2$. Moreover, $\lrb K_{k+1}\rrb$ is obviously constant along $T_{k+1}$ and is thus Dupin. On the other hand, as in the case of hypersurfaces we have that the point
\beeqs
\lrb K\rrb=\lrb\kappa Y_1+Y_{n+p+1}\rrb,\ \kappa\neq\infty,
\eneqs
is a curvature sphere of $\lambda$ at $\lp x,\xi\rp$ with corresponding principal vector $\lp X,-f_*A_\xi X+\eta\rp\in T_{\lp x,\xi\rp}UN_f$ if and only if
\beeqs
\begin{split}
\lp\kappa{Y_1}_*+{Y_{n+p+1}}_*\rp\lp X,-f_*A_\xi X+\eta\rp=\kappa\lp0,f_*X,0\rp+\lp0,-f_*A_\xi X+\eta,0\rp\\
=\lp0,-f_*\lp AX-\kappa X\rp+\eta,0\rp=0,
\end{split}
\eneqs
which implies that $\eta=0$ and that $\kappa$ is a principal curvature of $f$ with corresponding principal vector $X$, as we wished. Furthermore, we now have that
\beeqs
\begin{split}
K_*\lp X,-f_*A_\xi X\rp=\kappa_*\lp X,-f_*A_\xi X\rp Y_1+\lp\kappa{Y_1}_*+{Y_{n+p+1}}_*\rp\lp X,-f_*A_\xi X\rp\\
=\kappa_*\lp X,-f_*A_\xi X\rp Y_1,
\end{split}
\eneqs
which again is collinear with $K$ if and only if $\kappa_*\lp X,-f_*A_\xi X\rp=0$.
\enpr
As an immediate consequence of the above two results, we get a bijective correspondence between $k$-Dupin submanifolds of codimension $p\geq2$ (respectively, $k$-Dupin hypersurfaces) in space forms and $\lp k+1\rp$-Dupin (respectively, $k$-Dupin) Legendre lifts.
\becor\label{lgth}
Let $f:M^n\to\Q^{n+p}_c$ be an isometric immersion and $\lambda$ its Legendre lift. If $p\geq2$ (respectively, $p=1$), then $f$ is $k$-umbilical if and only if $\lambda$ is $\lp k+1\rp$-umbilical (respectively, $k$-umbilical). Furthermore, $f$ is $k$-Dupin if and only if $\lambda$ is $\lp k+1\rp$-Dupin (respectively, $k$-Dupin).
\encor
Another interesting application of Theorems \ref{cecil1} and \ref{cecil2} is the invariance of $k$-Dupin submanifolds under conformal transformations of and between ambient space forms.
\becor\label{cies}
Let $f:M^n\to\bar{\Q}^{n+p}_c$ be a $k$-umbilical isometric immersion and let $\phi:\bar{\Q}^{n+p}_c\to\bar{\Q}^{n+p}_{\tilde{c}}$ be a conformal diffeomorphism. Then, the composition $\tilde{f}=\phi\circ f:M^n\to\bar{\Q}^{n+p}_{\tilde{c}}$ is also $k$-umbilical. Moreover, if $f$ is $k$-Dupin, so is $\tilde{f}$.
\encor
\bepr
Since the Legendre lifts of $f$ and $\tilde{f}$ are M\"{o}bius conjugated by $\phi$, this follows immediately from the preceding corollary together with the fact that Lie sphere transformations map curvature spheres to curvature spheres and also preserve the Dupin condition.
\enpr
We close this subsection with a quick review of the theory Dupin submanifolds in the context of Lie sphere geometry. Pinkall \cite{Pinkall1985} introduced three constructions for obtaining a Dupin hypersurface in $\R^{n+m}$ from a Dupin hypersurface in $\R^{n+1}$, namely, tubes, cylinders and submanifolds of revolution, which give birth to a $\lp k+1\rp$-Dupin hypersurface out of a lower-dimensional $k$-Dupin hypersurface. These standard constructions can also be formulated in the context of Legendre submanifolds and apply even to Legendre lifts of submanifolds of higher codimension, although in this case the number of distinct curvature spheres remains the same.
\bere
When Pinkall introduced his constructions, he also listed cones, which turn out to be locally Lie equivalent to tubes.
\enre
A Dupin submanifold obtained from a lower-dimensional Dupin submanifold via one of Pinkall's standard constructions is said to be \emph{reducible}. More generally, a Dupin submanifold which is locally Lie equivalent to such a Dupin submanifold is called \emph{reducible}. Pinkall \cite{Pinkall1985} formulated the following simple Lie sphere geometric criterion for reducibility.
\bethr[Pinkall \cite{Pinkall1985}]\label{pith}
A connected proper Dupin submanifold $\lambda:N^{d-1}\to\Lambda^{2d-1}$ is reducible if and only if there exists a curvature sphere $\lrb K\rrb$ of $\lambda$ that lies in a linear subspace of $\mathbb{P}^{d+2}$ of codimension at least two.
\enthr
\bere
In \cite{dajczer2005}, Dajczer, Florit and Tojeiro studied reducibility in the context of Riemannian geometry. They formulated a concept of weak reducibility for proper Dupin submanifolds that have a flat normal bundle, including proper Dupin hypersurfaces. For hypersurfaces, their definition can be formulated as follows. A proper Dupin hypersurface $f:M^n\to\Q^{n+1}_c$ is said to be \emph{weakly reducible} if, for some principal curvature $\kappa_i$ with corresponding principal space $E_i$, the orthogonal complement $E_i^\perp$ is integrable. Dajczer, Florit and Tojeiro show that if a proper Dupin hypersurface $f:M^n\to\Q^{n+1}_c$ is Lie equivalent to a $\lp k+1\rp$-Dupin hypersurface that is obtained via one of the standard constructions from a $k$-Dupin hypersurface, then $f$ is weakly reducible. Thus, reducible implies weakly reducible for such hypersurfaces.

However, one can show that the open subset $U$ of a tube over the Veronese surface regarded in $\Sp^5$ (rather than in $\Sp^4$) on which there are three principal curvatures at each point is reducible but not weakly reducible, because none of the orthogonal complements of the principal spaces is integrable. Of course, $U$ is not constructed from a 2-Dupin submanifold, but rather one 3-Dupin.
\enre
A proper Dupin hypersurface with two distinct curvature spheres of respective multiplicities $p$ ad $q$ is called a \emph{cyclide of Dupin of characteristic} $\lp p,q\rp$. These are the simplest Dupin submanifolds after the spheres, and they were first studied in $\R^3$ by Dupin \cite{dupin1822} in 1822. An example of a cyclide of Dupin of characteristic $\lp1,1\rp$ in $\R^3$ is a torus of revolution. The cyclides were studied by many prominent mathematicians in the nineteenth century, including Liouville \cite{liouville1847}, Cayley \cite{cayley1873}, and Maxwell \cite{maxwell1867}, whose paper contains stereoscopic figures of the various types of cyclides. The long history of the classical cyclides of Dupin is given in Lilienthal \cite{lilienthal}. (See also Banchoff \cite{banchoff1970}, Cecil \cite{cecil1976}, Klein \cite{MR0226476}, Darboux \cite{MR1324110}-\cite{MR1365962}, Blaschke \cite{MR0015247}, Eisenhart \cite{eisenhart1909treatise}, Hilbert and Cohn-Vossen \cite{MR0046650}, Fladt and Baur \cite{MR0430974}, and Cecil and Ryan \cite{MR590644}, \cite{MR781126}, for more on the classical cyclides.) For a consideration of the cyclides in the context of computer graphics, see Degen \cite{MR1307650}, Pratt \cite{MR1074611}-\cite{MR1318323}, Srinivas and Dutta \cite{srinivas1994}-\cite{srinivas1994ii}-\cite{srinivas1995}-\cite{srinivas1995ii}, and Schrott and Odehnal \cite{MR2267352}. For cyclides of Dupin in $\R^3$, it was known in the nineteenth century that every connected Dupin cyclide is M\"{o}bius equivalent to an open subset of a surface of revolution obtained by revolving a profile circle $\Sp^1\subset\R^2$ about an axis $\R^1\subset\R^2\subset\R^3$. The profile circle is allowed to intersect the axis, thereby introducing Euclidean singularities. However, the corresponding Legendre map into the space of contact elements in $\R^3$ is an immersion. Higher-dimensional cyclides of Dupin appeared in the study of isoparametric hypersurfaces in spheres. Cartan knew that an isoparametric hypersurface in a sphere with two curvature spheres must be a standard product of spheres,
\beeqs
\Sp^d_{c_1}\times\Sp^{n-d}_{c_2}\subset\Sp^{n+1}\subset\R^{d+1}\times\R^{n-d+1}=\R^{n+2},\quad c_1^{-1}+c_2^{-1}=1.
\eneqs
Cecil and Ryan \cite{Cecil1978ii} showed that a compact cyclide of Dupin $M^n$ embedded in $\Sp^{n+1}$ must be M\"{o}bius equivalent to a standard product of spheres. The proof, however, uses the compactness of $M^n$ in an essential way. Later, Pinkall \cite{Pinkall1985} used Lie sphere geometric techniques to show that any two cyclides of Dupin of the same characteristic are locally Lie equivalent, which leads to a local classification of the higher-dimensional cyclides of Dupin that is analogous to the classical result, namely, every cyclide of Dupin is locally conformally congruent to a rotational hypersurface with a sphere as profile.
\bere
The latter result also follows easily from the Moore-type decomposition theorem for conformal immersions of a warped product of Riemannian manifolds due to Tojeiro \cite{Tojeiro2007}.
\enre
The case of 3-Dupin hypersurfaces was studied in the paper of Cecil and Jensen \cite{Cecil1998}. Their result can be stated as follows, where the case $d=4$ is credited to Pinkall \cite{Pinkall1985ii}.
\bethr[Pinkall \cite{Pinkall1985ii}, Cecil-Jensen \cite{Cecil1998}]\label{thcj}
Let $\lambda:N^{d-1}\to\Lambda^{2d-1},\ d\geq4$, be a connected irreducible 3-Dupin submanifold. Then $d=3\cdot2^{\mu-1}+1$, where $\mu=1,2,3,4$, each principal curvature has the same multiplicity $2^{\mu-1}$ and $\lambda$ is Lie equivalent to the Legendre lift of a 3-isoparametric hypersurface in $\Sp^{3\cdot2^{\mu-1}+1}$.
\enthr
\bere\label{reiv}
In \cite{Cartan1939}, Cartan showed that a 3-isoparametric hypersurface in the sphere is always a tube ${\psi_{\mathbb{F}}}_t:UN_{\psi_{\mathbb{F}}}^{3\cdot2^{\mu-1}}\to\Sp^{3\cdot2^{\mu-1}+1}$ over a standard embedding $\psi_{\mathbb{F}}:\mathbb{FP}^2\to\Sp^{3\cdot2^{\mu-1}+1}$ of one of the projective planes $\mathbb{FP}^2$, for $\mathbb{F}=\R,\mathbb{C},\mathbb{H},\mathbb{O}$ and $\mu=1,2,3,4$, respectively. Since ${\psi_{\mathbb{F}}}_t$ is exactly the spherical projection of the parallel submanifold $P_{-t}\lambda$ of the Legendre lift $\lambda:UN_{\psi_{\mathbb{F}}}^{3\cdot2^{\mu-1}}\to\Lambda^{3\cdot2^\mu+1}$ of $\psi_{\mathbb{F}}$, we have that $P_{-t}\lambda$ is the actual Legendre lift of ${\psi_{\mathbb{F}}}_t$. In other words, Legendre lifts of 3-isoparametric hypersurfaces and those of their corresponding standard embeddings of the projective planes are parallel, and in particular Lie equivalent.
\enre
\subsection{Generalized cylinders and envelopes of hyperspheres}\label{1pes}
This subsection is devoted to discussing the concepts of generalized cylinders and envelopes of hyperspheres.

Let $g:L^k\to\Q^{n+p}_c$ be an isometric immersion with a parallel flat normal subbundle $\mathcal{V}$ of rank $n-k$. The \emph{generalized $\lp n-k\rp$-cylinder} over $g$ determined by $\mathcal{V}$ is the $n$-dimensional submanifold parametrized by means of the exponential map of $\Q^{n+p}_c$ as
\beeqs
\gamma\in\mathcal{V}\mapsto\exp^c_{g\lp\pi\lp\gamma\rp\rp}\lp\gamma\rp,
\eneqs
where $\pi:\mathcal{V}\to L^k$ is the projection. It was shown by Dajczer-Florit-Tojeiro \cite{dajczer2005} that generalized cylinders are the only submanifolds that carry a relative nullity distribution with integrable conullity.

Let $f:M^n\to\Q^{n+1}_c$ be a connected orientable hypersurface, we denote by $\eta$ the unit normal to $f$ which gives the orientation of $M^n$. Following Asperti-Dajczer \cite{dajczer2005}, $f$ is said to be a $k$\emph{-parameter envelope of hyperspheres} (briefly, $k$-PES) if it carries a nonzero Dupin principal curvature $\kappa$ of multiplicity $n-k$. Actually, in \cite{dajczer2005} the authors restrict themselves to the case $k\geq2$, in which the Dupin condition is automatically satisfied.
Classically, a Euclidean $k$-PES is locally a solution
\beeqs
f=f\lp u_1,\dots,u_k,t_1,\dots,t_{n-k}\rp
\eneqs
in $\R^{n+1}$ of the system below:
\beeq\label{1pesdf}
\begin{aligned}
{} & \text{(a)}\quad\left\|f-g\right\|^2=r^2,\\
{} & \text{(b)}\quad\la f-g,\frac{\partial g}{\partial u_i}\ra=-r\frac{\partial r}{\partial u_i},\ i=1,\dots,k,
\end{aligned}
\eneq
where $g:L^k\to\R^{n+1}$ is an isometric immersion, $g=g\lp u_1,\dots,u_k\rp$, and $r\in C^\infty\lp L\rp$ is a non-vanishing function. Geometrically it means that $f$ is the envelope of the $k$-parameter family of hyperspheres given by (a): the limit of the intersection of neighboring hyperspheres that approach each other are $\lp n-k\rp$-spheres that generate the envelope.

Let $g$ and $r$ be as above and let $t_1,\dots,t_{n-k}$ be parameters of the unit $\lp n-k\rp$-sphere centered at the origin of a Euclidean $\lp n-k+1\rp$-space. Set
\beeq\label{1peseq}
f\lp u_1,\dots,u_k,t_1,\dots,t_{n-k}\rp=g-r\nabla r-r\sqrt{1-\left|\nabla r\right|^2}\varphi\lp t_1,\dots,t_{n-k}\rp,
\eneq
where $\nabla r$ is the gradient of $r$ and the vector $\varphi$ has origin at the point $\gamma=g-r\nabla r$ and describes a unit sphere in the affine $\lp n-k+1\rp$-plane through $\gamma$ orthogonal to $g$. It was shown by Asperti-Dajczer \cite{MR761994} that the hypersurface given by \eqref{1peseq} satisfies the system \eqref{1pesdf} and is (away from singular points) a $k$-PES. Conversely, every Euclidean $k$-PES satisfies system \eqref{1pesdf} and is locally of the form \eqref{1peseq}, for $r=1/\kappa$. Although their proof is carried out in Euclidean space, it can be easily adapted in nonflat space forms in order to get a similar parametrization of $k$-PES also in these ambient spaces.

Envelopes of hyperspheres have the following natural generalization in higher codimension. For simplicity, we discuss just the Euclidean case, the others being analogous. Given smooth maps $g:M^n\to\R^{n+p}$ and $r\in C^\infty\lp M\rp$, the family of hyperspheres
\beeqs
y\in M^n\mapsto S\lp g\lp y\rp;r\lp y\rp\rp
\eneqs
centered at $g\lp y\rp\in\R^{n+p}$ with radius $r\lp y\rp$ is said to be a \emph{congruence} of hyperspheres. It is said to be a $k$\emph{-parameter} congruence of hyperspheres if $g$ has rank $k$ everywhere and $\ker g_*\lp y\rp\subset\ker r_*\lp y\rp$ for all $y\in M^n$. An immersion $f:M^n\to\R^{n+p}$ is said to \emph{envelope} the congruence of hyperspheres along $M^n$ determined by $g:M^n\to\R^{n+p}$ and $r\in C^\infty\lp M\rp$ if
\beeqs
f_*T_yM\subset T_{f\lp y\rp}S\lp g\lp y\rp;r\lp y\rp\rp
\eneqs
for all $y\in M^n$.

A principal normal vector field $\eta$ is said to be \emph{Dupin} if $\eta$ is parallel in the normal connection along $E_\eta$. One can prove that if $f:M^n\to\R^{n+p}$ is an isometric immersion that carries a Dupin principal normal vector field $\eta$ of multiplicity $n-k$ then its focal submanifold
\beeq
g=f+\frac{1}{\left\|\eta\right\|^2}\eta
\eneq
and its curvature radius $r=\left\|\eta\right\|^{-1}$ determine an $\lp n-k\rp$-parameter congruence of hyperspheres that is enveloped by $f$. Conversely, every isometric immersion $f:M^n\to\R^{n+p}$  that envelopes an $\lp n-k\rp$-parameter congruence of hyperspheres along $M^n$ determined by $g:M^n\to\R^{n+p}$ and $r\in C^\infty\lp M\rp$ carries a Dupin principal normal vector field
\beeqs
\eta=\frac{1}{r^2}\lp g-f\rp
\eneqs
of multiplicity $k$.
\bere
Let $f:M^n\to\R^{n+1}$ be a $k$-PES. We say that $f$ is a \emph{special $k$-parameter envelope of hyperspheres} (briefly, $k$-SPES) if its Dupin curvature of multiplicity $n-k$ has integrable conullity. In \cite[Theorem (1.7)]{MR761994}, Asperti and Dajczer stated the following result on $k$-SPES:

\vspace{4pt}
\noindent\emph{ Let $f:M^n\to\R^{n+1}$ be a $k$-PES. Then $f$ is a $k$-SPES if and only if $g:L^k\to\R^{n+1}$ has flat normal bundle and $\nabla r\lp x\rp\in\Delta\lp x\rp$ for all $x\in L^k$, where $r=1/\kappa$.}

\vspace{4pt}
\noindent However, we point out that this result is not true. The counterexamples are the compact embedded cyclides of Dupin $f:M^n \to \mathbb{R}^{n+1}$ 
of characteristic $\lp p,q \rp$. In this case $f$ is both a $p$-SPES and a $q$-SPES. Since $M^n$ is compact then also $L^k$ is compact and thus the relative nullity $\Delta$ of the focal submanifold $g$ is trivial because $g$ is real analytic (this is so because $f$ is real analytic due to \cite{Cecil1978ii}). If the above statement were true the second condition would imply that these cyclides of Dupin are always isoparametric.

The mistake happens in the first step of their proof. The authors consider the general parametrization for $f$ given in equation \eqref{1peseq}. If in addition $f$ is a $k$-SPES, the authors claim that the local coordinates
\beeq\label{lcerr}
\lp u_1,\dots,u_k,t_1,\dots,t_{n-k}\rp
\eneq
in parametrization \eqref{1peseq} can be chosen such that $\lp u_1,\dots,u_k\rp$ are coordinates for the leaves of the Dupin conullity $E_\kappa^\perp$. Since $E_\kappa$ and $E_\kappa^\perp$ are integrable, there do exist local coordinates \eqref{lcerr} such that $\lp t_1,\dots,t_{n-k}\rp$ and $\lp u_1,\dots,u_k\rp$ are coordinates for the leaves of $E_\kappa$ and $E_\kappa^\perp$, respectively, but we cannot pick them without breaking the independence of $\varphi$ from the variables $u_1,\dots,u_k$. In other words, if the local coordinates \eqref{lcerr} are chosen as above, then parametrization \eqref{1peseq} becomes
\beeqs
f\lp u_1,\dots,u_k,t_1,\dots,t_{n-k}\rp=g-r\nabla r-r\sqrt{1-\left|\nabla r\right|^2}\varphi\lp u_1,\dots,u_k,t_1,\dots,t_{n-k}\rp.
\eneqs
\enre
\section{Proofs}\label{proofs}
In this section we prove our four main results stated in the introduction of this article. We start with the proof of Theorem \ref{2ith}, which is a simple application of Cartan's classification of 3-isoparametric hypersurfaces and the following result.
\bethr[Heintze-Olmos-Thorbergsson \cite{doi:10.1142/S0129167X91000107}]\label{cpcfi}
A submanifold of a space form has constant principal curvatures if and only if it is either isoparametric or a focal manifold to an isoparametric submanifold.
\enthr
For a discussion of focal manifolds of an isoparametric submanifold, its related Coxeter group, etc., see \cite{MR3468790}.
\bepr[Proof of Theorem \ref{2ith}]
Since $f$ has constant principal curvatures, as a consequence of Theorem \ref{cpcfi}, $f$ is either isoparametric or a focal manifold of an isoparametric submanifold. If $f$ is isoparametric, then it has flat normal bundle, which contradicts the assumptions that $f$ is 2-CPC and $p\geq2$. Thus, $f$ is a focal submanifold of an isoparametric submanifold $g:N^{n+q}\to\Q^{n+p}_c$. Then from the proof of Theorem \ref{cpcfi} it follows that $g$ is a holonomy tube $f_{\xi_x}$ through a principal vector $\xi_x\in N_fM\lp x\rp$. We claim that $g$ is an isoparametric hypersurface of a geodesic hypersphere. In fact, from \cite{MR3468790} it is enough to show that the codimension of a principal orbit $\mathcal{H}ol_x\lp\nabla^\perp\rp\cdot\xi_x$ is one. Indeed, if it is greater than one, we will obtain a contradiction. Note first that the Ricci equation implies $\lrb A_\xi,A_\eta\rrb=0$ for any two normal vectors $\xi,\eta$ to a principal orbit of the action of the normal holonomy group. So, we can perform a simultaneous diagonalization of the family $\lb A_\xi:\xi\in\text{the normal space to a principal orbit}\rb$. Then, given any two linearly independent normal vectors to a principal orbit, it is possible to find a nonzero umbilical linear combination of them, contradicting the assumption that $f$ is 2-CPC.

Thus, $g$ is an isoparametric hypersurface of a geodesic hypersphere and the isometric immersion $f$ is a focal submanifold of $g$, i.e., $f=g_{\xi/\kappa}$, where $\xi$ is the unit normal of $g$ and $\kappa$ a constant principal curvature of $g$. Since the shape operator $A_\xi$ of $f$ has two eigenspaces for any normal vector $\xi$ then the ``Tube formula'' (see \cite{MR3468790}) implies that the shape operator of $g$ must have three principal curvatures.

So from a theorem of Elie Cartan \cite{Cartan1939} it follows that $f$ must be one of the cited embeddings (for a beautiful proof of this Cartan's theorem see \cite{Console1998}).

Conversely, it is well known that the standard embeddings of the projective planes are 2-CPC.
\enpr
We now go on to prove our main Theorem \ref{omth}, for which we need three more lemmas. The first one asserts that the Legendre lift of a 2-Dupin submanifold of codimension greater than one is irreducible in the sense of Lie sphere geometry.
\bele\label{idl}
Let $f:M^{2m}\to\Sp^{2m+p},\ p\geq2,$ be a connected 2-Dupin isometric immersion. Then, its Legendre lift $\lambda:UN_f^{2m+p-1}\to\Lambda^{2\lp2m+p\rp-1}$ is an irreducible 3-Dupin submanifold with two principal curvatures of multiplicity $m$ and the other $p-1$.
\enle
\bepr
The fact that $\lambda$ is 3-Dupin follows directly from Theorem \ref{cecil2}. Moreover, the infinite principal curvature $\kappa_3$ has multiplicity $p-1$ and the two others have the same multiplicity $m$ by Lemma \ref{edl}.

It remains to be shown that $\lambda$ is irreducible. Let $\lrb K_i\rrb,\ 1\leq i\leq3,$ be the curvature spheres of $\lambda$ given by \eqref{cspc}, with $n=2m$. By Theorem \ref{pith} all there is left to prove is that no $\lrb K_i\rrb$ lies in a linear subspace of $\mathbb{P}^{2m+p+2}$ of codimension higher than one. For the point sphere map $\lrb K_3\rrb$, this follows easily from the fact that $f$ is conformally substantial. In fact, if $\lrb K_3\rrb$ lied in a $\mathbb{P}^{2m+p}$, then there would be a nonzero vector $v=\lp v_0,w,0\rp\in\R^{2m+p+3}_2$ other than $\lp0,0,1\rp$ that would be orthogonal to $\Span\lb K_3\rb$, so that
\beeqs
f\lp x\rp\cdot w=v_0,
\eneqs
which is the equation of a hyperplane in $\R^{2m+p+1}$. Thus, $f\lp M\rp$ would be contained in a hypersphere $\Sp^{2m+p-1}_c$, contradicting the fact that $f$ is conformally substantial.

To rule out the other two cases, note that, by Lemma \ref{edl}, $\kappa_i\lp x,-\xi\rp=-\kappa_j\lp x,\xi\rp,\ 1\leq i\neq j\leq2$, and thus
\begin{gather*}
\lrb K_i\lp x,\xi\rp\rrb=\lrb\kappa_iY_1\lp x,\xi\rp+Y_{2m+p+3}\lp x,\xi\rp\rrb=\lrb\lp\kappa_i,\kappa_if\lp x\rp+\xi,1\rp\rrb,\\
\lrb K_i\lp x,-\xi\rp\rrb=\lrb-\kappa_jY_1\lp x,-\xi\rp+Y_{2m+p+3}\lp x,-\xi\rp\rrb=\lrb\lp-\kappa_j,-\kappa_jf\lp x\rp-\xi,1\rp\rrb,
\end{gather*}
with $\kappa_i=\kappa_i\lp x,\xi\rp$ and similarly for $\kappa_j$. Adding these two equations, we get that
\beeqs
\begin{split}
\lrb K_i\lp x,\xi\rp+K_i\lp x,-\xi\rp\rrb=\lrb\lp\kappa_i-\kappa_j,\lp\kappa_i-\kappa_j\rp f\lp x\rp,2\rp\rrb\\
 =\lrb\lp\kappa_1-\kappa_2,\lp\kappa_1-\kappa_2\rp f\lp x\rp,\lp-1\rp^{i+1}2\rp\rrb\in\Span\lb K_{t_i}\rb,
\end{split}
\eneqs
for all $x\in M^{2m},\ 1\leq i\leq2$. Then the same argument used to discard $\lrb K_3\rrb$ gives that $f$ would reduce conformal codimension again if some $\lrb K_i\rrb$ lied in a linear subspace $\mathbb{P}^{2m+p}$, and this concludes the proof.
\enpr
\bere\label{ilex}
In the proof of Theorem \ref{w2uth} we will need a slightly more general version of Lemma \ref{idl}. Suppose that we are given a connected umbilic-free weakly 2-umbilical isometric immersion $f:M^{2m}\to\Sp^{2m+p},\ p\geq2$. Then, restricting the Legendre lift $\lambda$ of $f$ to the ``umbilic free'' part $UN_f^0$ of $UN_f$ (i.e., the set of $\lp x,\xi\rp\in UN_f$ such that $A_\xi$ is not a multiple of $I$), we still have by Theorem \ref{cecil2} and Remark \ref{le2ex}-(\emph{ii}) that $\lambda|_{UN_f^0}$ is 3-Dupin with the infinite principal curvature $\kappa_3$ of multiplicity $p-1$ and the two others of the same multiplicity $m$. Moreover, if we have in addition that $f$ is conformally substantial, then an entirely analogous proof shows that $\lambda|_{UN_f^0}$ is also irreducible.
\enre
The next two lemmas are related. Lemma \ref{btle} somehow characterizes the standard embeddings of the four projective planes as the only 2-umbilical submanifolds with codimension greater than one in space forms which can be focalized from a 3-Dupin hypersurface in two different ways. Lemma \ref{pl} is the Lie sphere geometric counterpart of Lemma \ref{btle} in terms of Legendre lifts.
\bele\label{btle}
Let $f:M^{2m}\to\Q_c^{2m+p},\ p\geq2,$ be a 2-umbilical isometric immersion. If the tube $f_t:UN_f^{2m+p-1}\to\Q_c^{2m+p}$ is nowhere an immersion for some $t\neq k\pi,\ k\in\mathbb{Z}$, if $c>0$ or $t\neq0$ if $c\leq0$, then actually $c>0$,
\beeqs
m=p-1=2^{\mu-1},\ \mu=1,2,3,4,
\eneqs
and, up to congruence and the similarity $\theta_c$, we have that
\beeqs
f\lp M\rp\subset\psi_{\mathbb{F}}\lp\mathbb{FP}^2\rp,
\eneqs
where $\mathbb{F}=\R,\mathbb{C},\mathbb{H},\mathbb{O}$, according to whether $\mu=1,2,3,4$, respectively.
\enle
\bepr
We will basically show that the assumption implies that $f$ is 2-unipotent, from which the statement will then follow, by Theorem \ref{dsth}.

Since the computations in all three cases are similar, we restrict attention to the spherical one ($c=1$). Then $f_t$ is given by
\beeq\label{tf}
f_t=\cos t\,f+\sin t\,\xi.
\eneq
The singularities of $f_t$ take place exactly at the focal points of $f$. Indeed, differentiating \eqref{tf}, we obtain
\beeq\label{tfd}
{f_t}_*\lp X,-f_*A_\xi X+\eta\rp=f_*\lp\cos t\,X-\sin t\,A_\xi X\rp+\sin t\,\eta,
\eneq
for all $X\in T_xM$ and all $\eta\in N_fM\lp x\rp$ orthogonal to $\xi$. Choosing $X=0$ in \eqref{tfd}, we get
\beeq
{f_t}_*\lp0,\eta\rp=\sin t\,\eta,
\eneq
but since $t\neq k\pi,\ k\in\mathbb{Z}$, it follows that $\sin t\neq0$, and thus we conclude that every $\eta\perp\xi$ lies in the range of ${f_t}_*$ at $\lp x,\xi\rp$. On the other hand, taking now $\eta=0$ in \eqref{tfd}, we have
\beeq\label{tfd2}
{f_t}_*\lp X,-f_*A_\xi X\rp=f_*\lp\cos t\,X-\sin t\,A_\xi X\rp,
\eneq
and so the assumption that $f_t$ is nowhere an immersion, together with the previous conclusion, forces \eqref{tfd2} to vanish for some $X\lp x,\xi\rp$ at every $\lp x,\xi\rp\in UN_f^{2m+p-1}$. In particular, it follows that $\cot t$ is a common principal curvature of all $A_\xi,\ \lp x,\xi\rp\in UN_f^{2m+p-1}$. Hence, we have that $f$ is 2-unipotent by Lemma \ref{upshl}. The result then follows immediately from Theorem \ref{dsth}.
\enpr
\bele\label{pl}
Let $f:M^{2m}\to\Q_c^{2m+p}, p\geq2,$ be a 2-umbilical isometric immersion and let $\lambda:UN_f^{2m+p-1}\to\Lambda^{2\lp2m+p\rp-1}$ be its Legendre lift. If the spherical projection of $P_t\lambda$ is nowhere an immersion for some spherical, Euclidean or hyperbolic parallel transformation $P_t\neq\pm I$, then
\beeqs
m=p-1=2^{\mu-1},\ \mu=1,2,3,4,
\eneqs
and, up to congruence and the stereographic projection $\pi_c$ if $c\leq0$ or similarity $\theta_c$ if $c>0$, we have that
\beeqs
f\lp M\rp\subset\psi_{\mathbb{F}}\lp\mathbb{FP}^2\rp,
\eneqs
where $\mathbb{F}=\R,\mathbb{C},\mathbb{H},\mathbb{O}$, according to whether $\mu=1,2,3,4$, respectively.
\enle
\bepr
After possibly composing $f$ with the stereographic projection, we are free to work in the ambient space to which the action of $P_t$ is best suited. More precisely, we can assume that $c=1,0,-1$ depending on whether $P_t$ is spherical, Euclidean or hyperbolic, respectively. Then the spherical projection of $P_t\lambda$ is simply the tube $f_{-t}$, and the result follows from the previous lemma.
\enpr

\subsection{Proof of Theorem \ref{omth}}

To prove the theorem all one has to do is put together the previous lemmata.

\bepr
By Lemma \ref{edl}, we have that $n=2m$. It follows from Lemma \ref{idl} that the Legendre lift $\lambda:UN_f^{2m+p-1}\to\Lambda^{2\lp2m+p\rp-1}$ of $f$ is a connected irreducible 3-Dupin submanifold with two principal curvatures of multiplicity $m$ and the other $p-1$. Thus, we are under the assumptions of Theorem \ref{thcj} and therefore we conclude that $p-1=m=2^{\mu-1}$ for $\mu=1,2,3,4$ (hence $2m+p=3\cdot2^{\mu-1}+1=4,7,13,25$, respectively) and $\lambda$ is Lie equivalent to an open subset of the Legendre lift of a 3-isoparametric hypersurface in $\Sp^4,\Sp^7,\Sp^{13},\Sp^{25}$, or equivalently, of that of the corresponding standard embedding $\psi_{\mathbb{F}}$ of the projective plane $\mathbb{FP}^2$, for $\mathbb{F}=\R,\mathbb{C},\mathbb{H},\mathbb{O}$, respectively (see Remark \ref{reiv}). In other words, there exists a Lie sphere transformation $\gamma$ such that the spherical projection $\tilde{f}$ of
\beeq\label{le}
\tilde{\lambda}=\gamma\lambda
\eneq
is a submersion of rank $n=2^\mu$ onto an open subset of the image of $\psi_{\mathbb{F}}$. On the other hand, we know from Theorem \ref{thcc} that any such $\gamma$ can be written as
\beeqs
\gamma=\Phi_1P_t\Phi_2,
\eneqs
where $\Phi_1$ and $\Phi_2$ are M\"{o}bius transformations and $P_t$ is some spherical, Euclidean or hyperbolic parallel transformation. Hence, taking $\Phi_1^{-1}$ on both sides of \eqref{le} we have
\beeqs
\Phi_1^{-1}\tilde{\lambda}=P_t\Phi_2\lambda.
\eneqs
Now let $\phi_i$ be the conformal transformation of $\Sp^{3\cdot2^{\mu-1}+1}$ induced by $\Phi_i$, $i=1,2$. Since $\lambda$ is the Legendre lift of $f$, it follows that
\beeqs
\Phi_2\lambda=\bar{\lambda}\bar{\phi}_2,
\eneqs
where $\bar{\phi}_2:UN_f^{3\cdot2^{\mu-1}}\to UN_{\phi_2\circ f}^{3\cdot2^{\mu-1}}$ is the unit normal bundle isometry induced by $\phi_2$ and $\bar{\lambda}:UN_{\phi_2\circ f}^{3\cdot2^{\mu-1}}\to\Lambda^{3\cdot2^\mu+1}$ is the Legendre lift of the composition $\phi_2\circ f$. Moreover, the spherical projection of $\Phi_1^{-1}\tilde{\lambda}$ is simply $\phi_1^{-1}\circ\tilde{f}$, and thus is also a submersion of rank $2^\mu$ onto its image. Hence, so must be that of $P_t\bar{\lambda}$ and, in particular, is nowhere an immersion, from which it follows, by Lemma \ref{pl}, that already $\phi_2\circ f\lp M\rp\subset\psi_{\mathbb{F}}\lp\mathbb{FP}^2\rp$ up to congruence, or else $P_t=\pm I$. In both cases we conclude that
\beeq\label{cnt}
f\lp M\rp\subset\psi_{\mathbb{F}}\lp\mathbb{FP}^2\rp
\eneq
after a conformal motion, which completes the proof of the direct statement.

The converse follows immediately from the conformal invariance of the class of 2-Dupin submanifolds given by Corollary \ref{cies}.
\enpr
\subsection{Proofs of Theorems \ref{2uth} and \ref{w2uth}}
Here we prove Theorem \ref{2uth}.
\bepr[Proof of Theorem \ref{2uth}]
If $n=2$, then $p\leq2$, otherwise there are umbilical directions, since the vector subspace of real symmetric $2\times2$ matrices is 3-dimensional. This is case (\emph{iv}) in the statement. Now we can assume $n\geq3$. If every principal curvature has multiplicity greater than or equal to two, then $f$ is automatically 2-Dupin, so cases (\emph{i}) and (\emph{ii}) follow from Theorem \ref{omth} and Pinkall's classification of the cyclides, respectively. Now we can assume that one of the principal curvatures has multiplicity one. From Lemma \ref{le2ex} it follows that $p=1$. Let $\kappa$ be the principal curvature of multiplicity $n-1$. Then, on each connected component of the open subset where $\kappa\neq0$, we have by definition that $f$ is a 1-PES. Finally, on each connected component of the interior of the closed subset where $\kappa=0$, we have by the result of Dajczer-Florit-Tojeiro \cite{dajczer2005} that $f$ is a generalized $\lp n-1\rp$-cylinder (see Subsection \ref{1pes}).
\enpr
Now we prove Theorem \ref{w2uth}.
\bepr[Proof of Theorem \ref{w2uth}]
In the case $n=2$ there is nothing to be done. So we assume $n\geq3$. Assume that the open subset
\beeqs
V=\lb x\in M^n:R^\perp\lp x\rp\neq0\rb
\eneqs
is nonempty. By Remark \ref{le2ex}-(\emph{ii}), we have that $n=2m$, with $m\geq2$.

Take a connected component $U$ of $V$, let $q$ be the conformal codimension of $f|_U$ and regard $f|_U$ as an isometric immersion into a totally umbilical submanifold $\Sp^{2m+q}_c$ of $\Sp^{2m+p}$, with $c\geq1$. Then Remark \ref{ilex} implies that its umbilic-free Legendre lift $\lambda:UN_f^0\to\Lambda$ is a connected irreducible 3-Dupin submanifold with two curvature spheres of multiplicity $m$ and the other $q-1$. Thus, we are under the assumptions of Theorem \ref{thcj} and therefore we conclude that $q-1=m=2^{\mu-1}$ for $\mu=2,3,4$ (hence $2m+q=3\cdot2^{\mu-1}+1=7,13,25$, respectively) and $\lambda$ is Lie equivalent (by a Lie sphere transformation $\gamma$) to an open subset of the Legendre lift of the standard embedding $\psi_{\mathbb{F}}$ of the projective plane $\mathbb{FP}^2$, for $\mathbb{F}=\mathbb{C},\mathbb{H},\mathbb{O}$, respectively.  As in the proof of Theorem \ref{omth}, by decomposing $\gamma$ according to Theorem \ref{thcc}, we get that
\beeq\label{feqf}
f\lp U\rp\subset\text{Im}\,\psi_{\mathbb{F}}
\eneq
after a conformal motion.

Finally, we claim that $U=M$. Assume otherwise and take $x\in\partial U$. Since $\psi_{\mathbb{F}}$ is homogeneous and has nonflat normal bundle, it follows that $\left\|R_{\psi_{\mathbb{F}}}^\perp\right\|=\kappa$ is constant and positive, but then \eqref{feqf} yields
\beeqs
\left\|R_f^\perp\lp x\rp\right\|=\kappa>0,
\eneqs
which is a contradiction with $R_f^\perp\lp x\rp=0$. Thus the claim is proved and so
\beeq
\text{Im}\,f\subset\text{Im}\,\psi_{\mathbb{F}},
\eneq
which gives us case (\emph{i}) in the statement.

Now we can assume that $R^\perp=0$ everywhere. If the two principal normals $\eta_1$ and $\eta_2$ of $f$ have multiplicity $\geq2$, then they are automatically Dupin, and hence the orthogonal net $\mathcal{E}=\lp E_1,E_2\rp$ of their principal distributions is what has been called a CWP-net in \cite{Tojeiro2007}, which implies that $M^n$ is locally conformal to a warped product. Moreover, since the second fundamental form $\alpha$ of $f$ is adapted to $\mathcal{E}$ (i.e., $\alpha\lp E_1,E_2\rp=0$) and also the leaves of $E_i,\ 1\leq i\leq2,$ are extrinsic spheres of $\Sp^{n+p}$, case (\emph{ii}) then follows directly from the main theorem in \cite{Tojeiro2007}.

Finally, case (\emph{iii}) i.e. when one of the multiplicities is one, is obtained from the considerations of Subsection \ref{1pes} in a similar way as in the proof of Theorem \ref{2uth}.
\enpr
\bibliographystyle{utphys}
\bibliography{ourbibumb}
\end{document}